\documentclass[a4paper,11pt]{amsart}
\usepackage[T1]{fontenc}
\usepackage[utf8]{inputenc}
\usepackage[foot]{amsaddr}

\usepackage[pdftex]{hyperref}
\usepackage[margin=2.5cm]{geometry}
\usepackage[english]{babel}
\usepackage{amsmath,amssymb,verbatim,mathrsfs,latexsym,paralist}
\usepackage[dvipsnames]{xcolor}
\usepackage{graphicx}
\usepackage{subcaption}
\usepackage{epstopdf}
\usepackage{indentfirst}
 \usepackage{amsthm}
 \allowdisplaybreaks[4]
 \usepackage{tikz}
\usepackage{caption}
\usepackage[T1]{fontenc}
\usepackage[all]{xy}
\usepackage{xcolor}

\newcommand{\Z}{\mathbb{Z}}
\newcommand{\C}{\mathbb{C}}

\newcommand{\Q}{\mathbb{Q}}

\newcommand{\mrm}{\mathrm}
\newcommand{\mf}{\mathfrak}

\newcommand{\pos}{\mrm{pos}}
\newcommand{\mut}{\mrm{Mut}}

\newcommand{\Ext}{\mrm{Ext}}
\newcommand{\Hom}{\mrm{Hom}}

\numberwithin{equation}{section}
\newtheorem{theorem}{Theorem}[section]
\newtheorem*{theorem*}{Theorem}
\newtheorem*{thm}{Theorem}
\newtheorem{proposition}[theorem]{Proposition}

\newtheorem{lemma}[theorem]{Lemma}
\newtheorem{corollary}[theorem]{Corollary}
 
\newtheorem{defn}[theorem]{Definition}  
\newtheorem{example}[theorem]{Example}   
\newtheorem{crit}[theorem]{Criterion}
\newtheorem*{lemma*}{Lemma}
\newtheorem*{defn*}{Definition} 

\theoremstyle{remark}
\newtheorem{rmk}[theorem]{Remark}

\title{On the Smooth Locus in Flat Linear Degenerations of Flag Varieties}
\author[Sabino Di Trani]{Sabino Di Trani}
\address{ Dipartimento di Matematica ``Guido Castelnuovo'', Sapienza - Universit\`a di Roma.
}
\email{sabino.ditrani@uniroma1.it }
\address{\emph{The author has been partially supported by GNSAGA - INDAM group.}}
\address{ORCID id: https://orcid.org/0000-0002-6651-558X}
\begin{document}

\maketitle
\textbf{Abstract:}
We use some moment graph techniques, recently introduced by Lanini and Pütz, to provide a combinatorial description of the smooth locus in flat linear degenerations of flag varieties, generalizing a result proved by Cerulli Irelli, Feigin and Reineke for the Feigin
degeneration. Moreover, we propose a different combinatorial criterion linking the smoothness at a fixed point to transitive tournaments. \\

\section{Introduction}
 The \emph{Complete Flag Variety} $\mathcal{F}l(\C^{n+1})$ is a very well known object, playing a relevant role in many different areas of mathematics. It is a  smooth irreducible algebraic variety of dimension $\binom{n+1}{2}$ and can be naturally realized as  the variety of $n$-tuples of nested subspaces  $(V_1, \dots , V_n)$ of $\C^{n+1}$ such that $\dim V_i = i$.\\
From a slightly different point of view, each subspace $V_i$ of an $n$-tuple $(V_1, \dots, V_n)$ can be thought as a subspace of its own copy $W_i$ of $\C^{n+1}$ and consequently it is possible to identify the complete flag variety with the closed subvariety of points $(V_1, \dots , V_n) $ in the product of Grassmann varieties $\mrm{Gr}_1(W_1)\times \dots \times \mrm{Gr}_n(W_n)$, that satisfy the relations $\mrm{id}_i V_i \subset V_{i+1}$ for every~$i$, where $\mrm{id}_i$ denotes the identity map from $W_i$ to $W_{i+1}$.\\
It is possible to generalize this construction 
considering any $n-1$-tuple $\mathbf{f}:=(f_1, \dots, f_{n-1})$ of linear endomorphisms of $\C^{n+1}$ and prescribing that $f_i V_i \subset V_{i+1}$.
\begin{defn*}
For a fixed family $\mathbf{f}=(f_1, \dots, f_{n-1})$ of linear endomorphisms of $\C^{n+1}$, the \emph{Linear Degeneration  of Flag Variety associated to the family of maps $\mathbf{f}$} is the algebraic variety of $n$-tuples of subspaces $(V_1, \dots , V_n)$ such that $\dim V_i = i$ and $f_i V_i\subset V_{i+1}$. We denote it by~${\mathcal{F}l}^\mathbf{f}(\C^{n+1})$. 
\end{defn*}
We also refer to ${\mathcal{F}l}^\mathbf{f}(\C^{n+1})$ as an $\mathbf{f}$-degeneration of the complete flag variety.
The ubiquitous appearing of linear degenerations in many different context of mathematics as representation theory, algebraic geometry and combinatorics, motivated a growing interest for these objects, that are extensively studied in a large amount recent papers (see  \cite{Fe}, for an introductory survey about linear degenerations and their role in mathematics, and \cite{CIII}, for a complete and more technical exposition on this topic). \\
One of the crucial points in studying linear degenerations is that their appear as a special case of Quiver Grassmanians, for certain representations of the equioriented quiver of type $A_n$. More precisely, to the family of linear endomorphisms $\textbf{f}=(f_1, \dots, f_{n-1})$ it is possible to associate an $A_n$-representation considering
$$M_\textbf{f}:= \quad  \C^{n+1} \stackrel{{f_1}}{\longrightarrow} \C^{n+1} \stackrel{{f_2}}{\longrightarrow} \ldots  \stackrel{{f_{n-2}}}{\longrightarrow}\C^{n+1} \stackrel{{f_{n-1}}}{\longrightarrow}\C^{n+1}. $$ 
An $n$-tuple of subspaces $(V_1, \dots, V_{n})$ with $\mrm{dim}V_i=i$ and such that $f_i V_i \subset V_{i+1}$ corresponds to a subrepresentation of $M_\textbf{f}$ with dimension vector $\textbf{e}=(1, \dots, n)$ and consequently the linear degeneration associated to $\textbf{f}$ can be identified with the quiver Grassmanian $\mrm{Gr}_{\textbf{e}}(M_\textbf{f})$.
This interpretation of linear degenerations of flag variety from a "Quiver Representation Theory"- point of view has the merit of highlighting how relevant geometric properties of these objects can be recovered using a combination of both combinatorial tools and  results of algebraic geometry. \\
Using their realization as quiver Grassmanians, it is possible to endow linear degenerations with suitable torus actions,  allowing very often to obtain deep results about their topology and explicit computations concerning their homological invariants. \\
As an example, using these techniques in~\cite{CIIII} it is proved that the Euler characteristics of the Feigin degenerations, i.e. degenerations of ${\mathcal{F}l}(\C^{n+1})$ such that each $f_i$ has rank $n$ and the kernels of $f_1, \dots, f_{n-1}$ are linearly independents, are equal to the median Genocchi numbers. \\
Moreover, in~\cite{CII} a decomposition of Feigin degenerations into disjoint cells with respect to a suitable torus action is achieved and smoothness of points in a cell can be verified looking at certain sequences of integers (see \emph{CI-F-R Condition}, Definition \ref{CIFR}). \\
More recently, actions of some bigger tori on linear degenerations were introduced and studied in \cite{MA} and \cite{MAII}, as a special case of torus actions on Cyclic Quiver Grassmanians. \\
As a consequence of constructions contained in \cite{MA}, every $\mathbf{f}$-degeneration ${\mathcal{F}l}^\mathbf{f}(\C^{n+1})$ can be equipped with the action of a suitable torus  $T$, such that the $T$-action has a finite number of fixed points and a finite number of 1-dimensional orbits. It also results that ${\mathcal{F}l}^\mathbf{f}(\C^{n+1})$ has a Bialynicki-Birula decomposition into rational cells, induced by a generic $T$-cocharacter. The degeneration ${\mathcal{F}l}^\mathbf{f}(\C^{n+1})$ is then a \emph{GKM Variety}~(Definition \ref{defn:GKM}) and its $T$-equivariant cohomolo\-gy  can be determined using the \emph{Moment Graph}
~(Definition~\ref{defn:moment}). Moreover the moment graph encodes some relevant information about the dimension of tangent space at $T$-fixed points. \\
An explicit description of the moment graph in combinatorial terms is provided in \cite{MA} and \cite{MAII}. More precisely, every fixed point for the $T$ action on ${\mathcal{F}l}^\mathbf{f}(\C^{n+1})$ can be associated to a family $S=(S_1, \dots S_n)$ of subsets of \{1,\dots, n+1\}, called \emph{Admissible Sequence} (Definition~\ref{defn:admissible}) and edges in the moment graph correspond to combinatorial operations, called \emph{Mutations} (Definition~\ref{defn:mutation}), on certain quivers associated to $T$- fixed points.\\
Inspired by one of the many open questions proposed by Lanini and P\"{u}tz (\cite{MA}, Introduction, Question~(5)), the main purpose of this paper is to use GKM theory to provide smoothness criteria for $T$-fixed points in flat linear degenerations, i.e. linear degenerations that are equidimensional varieties of dimension~$\binom{n+1}{2}$. Some of these criteria appear as a natural combinatorial generalization of the CI-F-R Condition to all flat linear degenerations. Moreover, our criteria highlights also a deeper link between smoothness and combinatorics of mutations.\\
\textbf{Organization of the Paper:} In Section~\ref{sec:2} linear degenerations of flag variety are presented in a more geometric flavour as fibers of a suitable $GL_{n+1}(\C)^n$-equivariant map from the \emph{Universal Degenerate Flag Variety} to the total space $\mrm{Hom}(\C^{n+1}, \C^{n+1})^{n-1}$ of representation of the equioriented quiver $A_n$. In this geometrical setting \emph{Flat Degenerations}~(Definition~\ref{defn:flat}) can be introduced in a very natural way and they can be classified completely up to $GL_{n+1}(\C)^n$-action. \\  
Section~\ref{sec:quivers} contains all the relevant informations about quiver representations that we are going to use to prove our results.
In particular we focus on the construction of \emph{Coefficient Quiver~$Q(M)$} (Definition~\ref{defn:coeffquiv}) associated to a representation $M$ of the equioriented quiver of type $A_n$. \\
In Section \ref{sec:3} we summarize the results proved in \cite{MA} and \cite{MAII} about the structure of GKM varieties of linear degenerations. In particular we recall how the structure of the moment graph can be described  using mutations between \emph{Successor Closed Subquivers}~(Definition~\ref{defn:scs}) of the coefficient quiver $Q(M_\textbf{f})$.\\
In Section~\ref{sec:formula} we prove a formula for the dimension of tangent space at a fixed point $p_S$, that involves certain combinatorial properties of the admissible sequence $S=(S_1, \dots S_n)$. More precisely, for each point $p_S$ we define a set $\mathrm{Sing}(p_S)$, that encodes in a combinatorial way certain homological informations about indecomposables appearing in $p_S$, when viewed as subrepresentation of $M_\textbf{f}$. In particular we prove that $|\mathrm{Sing}(p_S)|$ equals the dimension of $\mathrm{Ext}^1(p_S, M_\textbf{f}/p_S)$ and we obtain that\[\dim T_{p_S}{\mathcal{F}l}^\textbf{f}(\C^{n+1}) = \frac{n(n+1)}{2}+ |\mathrm{Sing} (p_S)|\] as a consequence of Euler Formula.\\
 Section \ref{sec:4} is devoted to presenting our smoothness criteria for $T$-fixed points in flat degenerations. We identify a combinatorial property of admissible sequences, the \emph{Generalized CI-F-R Condition}~(Definition~\ref{CIFR}), that generalizes the
 smoothness condition provided in \cite{CII} for Feigin degeneration. Moreover, to each point $p_S$ we attach an oriented graph, the \emph{Oriented Mutation Graph} $\widetilde{G}_S$ ~(Definition~\ref{defn:mutG}). This graph encodes relevant informations about 
vertices adjacent to $p_S$ in the moment graph. In our main result, we prove that smoothness of $p_S$ can be directly checked looking a these combinatorial objects:
\begin{theorem*}[Smoothness Criteria] Let $p_S$ be a $T$ fixed point in a flat linear degeneration of ${\mathcal{F}l}(\C^{n+1})$, associated to an admissible sequence $S$. The following conditions are equivalent:
\begin{enumerate}
    \item The point $p_S$ is smooth;
    \item The set $\mrm{Sing}(p_S)$ is empty;
        \item The admissible sequence $S$ has the Generalized CI-F-R Condition;
            \item The Oriented Mutation graph $\widetilde{G}_S$ is a transitive tournament $n+1$ vertices;
    \item The unoriented graph underlying $\widetilde{G}_S$ is the complete graph over $n+1$ vertices.
\end{enumerate}
\end{theorem*}
At the end of the section, some counterexamples to our criteria are provided, in the case of non flat degenerations.
Finally, Section~\ref{sect:proofs} contains the proofs of our results.\\
 \textbf{Acknowledgments:} I would like to express my deep gratitude to Martina Lanini, for introducing me to combinatorics of $T$-actions on quiver Grassmanians and for her constant support during the preparation of this work. Moreover, I am grateful to Giovanni Cerulli Irelli for his comments to a first version of this paper and for many useful discussions about linear degenerations and quiver representations. Finally, I have to thank Marco Trevisiol for some helpful remarks about the organization of the paper. 
\section{A Geometric Picture}\label{sec:2}
In this section we recall some well known geometrical properties of linear degenerations.
We refer to \cite{CIII} for a complete treatment on this subject. Consider the two algebraic varieties 
\[U := \mrm{Hom}(\C^{n+1}, \C^{n+1})^{n-1} \qquad Z:=\mrm{Gr}_1(\C^{n+1})\times \dots \times \mrm{Gr}_n(\C^{n+1}). \]
The set
\[Y=\{(\mathbf{f}, \mathbf{V})  \, | \; \mathbf{f}=(f_1, \dots, f_{n-1}) \in U, \, \mathbf{V}=(V_1, \dots, V_n) \in Z, \, f_i V_i \subset V_{i+1}\}\]
 is a smooth closed subvariety of $U \times Z$ and is known as the \emph{Universal Linear Degeneration of the Flag Variety}.  The group $G:=\mrm{GL}_{n+1}(\C)^n$ acts on $U$ by the rule 
\[(g_1, \dots, g_n)\cdot (f_1, \dots, f_{n-1}):=(g_2f_1g_1^{-1}, \dots, g_nf_{n-1}g_{n-1}^{-1}) \]
and on $Z$ componentwise, by the classical matrix action on $\C^{n+1}$. So $G$ acts on the product variety $U \times Z$ too, and the subvariety $Y$ can be endowed of a structure of $G$-variety in a way that the two projections $\pi_1: Y \rightarrow U$ and $\pi_2: Y \rightarrow Z$ result to be $G$-equivariant morphisms.\\ 
The space $U$ can be thought as a parameter space for $Y$: to any $n-1$-tuple of morphisms $\mathbf{f}=(f_1, \dots, f_{n-1})$, it is associated a closed subvariety of $Y$, considering the fibre $\pi_1^{-1}(\mathbf{f})$. For a fixed $\mathbf{f}$, this variety is isomorphic to the degeneration ${\mathcal{F}l}^\mathbf{f}(\C^{n+1})$.\\ Moreover, the $G$-equivariance of $\pi_1$ implies that the isomorphism class of each degeneration is uniquely determined by the sequence $\{r_{i,j}(\mathbf{f})\}$ of ranks of composition maps $f_i\circ \dots \circ f_j$. \\
Geometric properties of linear degenerations of the flag variety can be checked directly by looking at properties of $\pi_1$, extensively investigated in \cite{CIII}. In particular the authors prove that there exist two maximal subsets $U_{irr, \,flat} \subset U_{flat}\subset U$ such that 
\begin{itemize}
    \item the restriction of $\pi_1$ to $\pi_1^{-1} \left(U_{flat}\right)$ is a flat morphism of algebraic varieties;
    \item the fiber $\pi_1^{-1}(\mathbf{f})$ is an irreducible algebraic variety for every $\mathbf{f} \in U_{irr, \,flat}$
\end{itemize}
Complete flag variety~${\mathcal{F}l}(\C^{n+1})$ appear in $\pi_1^{-1} \left( U_{flat} \right)$ as the fiber of the element $\textbf{id}=(\mrm{id}, \dots, \mrm{id})$. Flatness then implies that for every $\mathbf{f} \in U_{flat}$ the irreducible components of  $\pi_1^{-1}(\mathbf{f})$ are equidimensional of dimension $n(n+1)/2$. Moreover, in \cite{CIII} it is proved that this property characterizes exactly the linear degenerations appearing as fibers over $U_{flat}$.
\begin{defn}\label{defn:flat}
An $\mathbf{f}$-degeneration ${\mathcal{F}l}^\mathbf{f}(\C^{n+1})$, where $\mathbf{f} \in  U_{flat}$ (resp. $\mathbf{f}\in  U_{irr, flat}$), is called a \emph{Flat} (resp. \emph{Flat Irreducible}) \emph{Degeneration} of the Flag Variety.
\end{defn}
Sets  $U_{irr, \,flat}$ and $U_{flat}$ are completely characterized in \cite{CIII} in terms of rank sequences~$\{r_{i,j}(\mathbf{f})\}$. 
\subsection{Classification of Flat Degenerations}
For each $G$ orbit in the flat locus, we want to identify a suitable orbit representative. \\
Let us fix a basis $B=\{v_1, \dots, v_{n+1}\}$ of $\C^{n+1} $ and  $R \subset \{1, \dots, n+1\}$, the projection operator $\pi_R: \C^{n+1} \rightarrow \C^{n+1}$ (with respect to $B$) is the linear operator 
on $\C^{n+1} $ 
defined by the rule: 
\[\pi_R(v_i)=\begin{cases} 0 & \text{if } i \in R,  \\ v_i & \text{ otherwise } \end{cases}\]
For any family $\mathbf{R}=(R_1, \dots, R_{n-1})$ of subsets of $\{1, \dots, n+1\}$ we define the degenerate flag variety ${\mathcal{F}l}^\mathbf{R}(\C^{n+1} )$ associated to $\mathbf{R}$ as the $\mathbf{f}$-degeneration of ${\mathcal{F}l}(\C^{n+1} )$ obtained considering the family of endomorphisms $\mathbf{f}=(f_1, \dots, f_{n-1})$ such that 
$f_i=\pi_{R_i}$ for every $i$.\\
Each $G$-orbit has a representative of the form ${\mathcal{F}l}^\mathbf{R}(\C^{n+1} )$, so in the remaining of the paper we reduce our analysis to studying these special degenerations.\\
Using some elementary linear algebra, in \cite{CIII} orbit representatives for flat degeneration are described in terms of certain sequences of integers. 
\begin{thm}[c.f.r.\cite{CIII}, Remark 2]\label{thm:flat}
Let $\mathbf{R}=(R_1, \dots, R_n)$ be a family of subsets of $\{1, \dots, n+1\}$. The degeneration ${\mathcal{F}l}^\mathbf{f}(\C^{n+1})$ is flat if and only if its orbit has a representative of the form ${\mathcal{F}l}^\mathbf{R}(\C^{n+1} )$, where $\mathbf{R}$ satisfies the following conditions for all $ h \in \{1, \dots, n\}$: 
\begin{enumerate}\item  $|R_h|\leq 2$; \item $|R_h\cup R_{h+1}| \leq 3$ \end{enumerate}
Moreover, ${\mathcal{F}l}^\mathbf{R}(\C^{n+1} )$ is a flat irreducible degeneration if and only if $|R_h|\leq 1$ for all $h$.
\end{thm}
 \begin{example}[The Feigin Degeneration] An $\mathbf{f}$-degeneration such that each $f_i$ has rank $n$ and such that the kernels of $f_1, \dots, f_{n-1}$ are linearly independents is a \emph{Feigin Degeneration}. Up to base change every Feigin degeneration is of the form ${\mathcal{F}l}^\mathbf{R}(\C^{n+1} )$, where $\mathbf{R}$ is the family of subsets~$(\{2\},\dots,\{n-1\})$. 
It is the toy model for our results. We summarize here some of its remarkable properties (see \cite{CIII}, \cite{CIIII} and \cite{CII} for some more complete and precise statements): 
\begin{itemize}
    \item Feigin Degeneration is a singular, normal, irreducible, local complete intersection, projective variety of dimension $n(n+1)/2$;
    \item It can be endowed of a structure of $T$-variety for a suitable algebraic torus $T$. The $T$ action has a finite number of fixed points and of 1-dimensional orbits;
    \item The torus action induces a cellular decomposition of the Feigin Degeneration into a finite number of disjoint cells. Each cell contains exactly one fixed point;
    \item It can be defined an action of an algebraic group $\mf{A}$ on the Feigin degeneration such that each cell coincides with the orbit of the associated $T$-fixed point under the $\mf{A}$-action.
\end{itemize}
\end{example}
Because of cellular decomposition of Feigin degeneration into orbits of fixed points, it is natural to look for a criterion to determine if a fixed point is smooth. Such a criterion is proved in~\cite{CII}. We will recall it in Section \ref{sec:4}, as a starting point for our more general results.
\section{Quiver Representations}\label{sec:quivers}
In this section we summarize some results about representations of a generic quiver $Q$. For a complete reference on this subject we refer to \cite{ASS} and \cite{Ki}.  Moreover we recall how representations of quivers of Type $A$ are related to linear degenerations of flag varieties. \\ Let $Q$ be a finite quiver with set of vertices $Q_0$ and set of oriented edges $Q_1$. If $\alpha \in Q_1$, we denote by $s(\alpha)$ and $t(\alpha)$ the source and the target of $\alpha$, respectively. 
\begin{defn}
A representation $(M,F)$ with base field  $\C$ of a quiver $Q=(Q_0,Q_1)$ is the datum of:
\begin{itemize}
    \item a family of finite dimensional complex vector spaces $M=(M_i)_{i \in Q_0}$,
    \item a sequence of maps $F=(f_\alpha)_{\alpha \in Q_1 }$ such that $f_\alpha \in \mrm{Hom}_\C(M_{s(\alpha)}, M_{t(\alpha)})$.
\end{itemize}
\end{defn}
Often we omit the family $F$ in the notation if the context is clear. To each representation $M$ is attached a dimension vector ${\bf{d}}=(d_i)_{i \in Q_0}$ such that $d_i=\mrm{dim}_\C M_i$.
A subrepresentation $N$ of $M$ is a collection of subspaces $N_i \subset M_i, \; i \in Q_0$ compatible with the maps~$f_\alpha$. 
\begin{defn}
Let $M$ be a representation of a quiver $Q$. Let ${\bf{e}}=(e_i)_{i \in Q_0}$ be a vector of non negative integers. The Quiver Grassmanian $\mrm{Gr}_{\bf e}(Q, M)$ is the algebraic variety of subrepresentations of $M$ with dimension vector $\bf e$. 
\end{defn}
\begin{example}\label{DegRap}Consider the quiver $Q$ such that $Q_0=\{1, \dots, n\}$ and $Q_1=\{(i,i+1) | 1 \leq i < n\}$. We will refer to this quiver as the equioriented quiver of type $A_n$. \begin{enumerate}\item
The assignment $M_i=\C^{n+1}$ for all $i \in Q_0$ and $\varphi_\alpha=\mathrm{id}$ for all $\alpha \in Q_1$ defines a representation $M$ of $A_n$. Each subrepresentation of $M$ of dimension vector $(1, \dots, n)$ corresponds to a point in ${\mathcal{F}l}(\C^{n+1})$. \item
Fix family of endomorphisms $\mathbf{f}=\{f_1, \dots , f_n\} $ of $\C^{n+1}$ and  consider the representation $M_{\bf{f}}$ as described in the introduction. 
If we set $\textbf{e}=(1, \dots, n)$, the quiver Grassmanian $\mrm{Gr}_{\bf e}(A_{n}, M_{\bf{f}})$ is isomorphic to the linear degeneration  ${\mathcal{F}l}^\mathbf{f}(\C^{n+1})$.\end{enumerate}
\end{example}
In our context we always work with representation of equioriented quivers of type $A$ so we omit the quiver $Q$ in the quiver Grassmanian notation.

We denote by $\mrm{Rep}(Q)$ the category of representations of the quiver $Q$ over $\C$. A general fact about quiver representations is that $\mrm{Rep}(Q)$ is hereditary, i.e. $\Ext^{\geq 2}(M,N)=0$ for every $M,N \in \mrm{Rep}(Q)$. As a consequence, if $M$ and $N$ are representations of dimension vectors $\bf d$ and~$\bf e$ respectively, the following formula holds:
\begin{equation}\label{ExtHom}
\dim_\C \Hom_Q(M,N)-\dim_\C \Ext^1(M,N)=\langle \bf d, \bf e \rangle,  
\end{equation}
where $\langle - \, , \, -  \rangle$ denotes the Euler Form associated to the quiver $Q$.
\begin{rmk}If $Q$ is the equioriented quiver of type $A_n$, the Euler Form is defined by the formula
$$\langle{\bf d},{\bf e}\rangle=\sum_{i=1}^nd_ie_i-\sum_{i=1}^{n-1}d_ie_{i+1}.$$ \end{rmk}
As shown in Example~\ref{DegRap}, each linear degeneration is isomorphic to a quiver Grassmanian, so it is natural to use techniques coming from representation theory of quivers of type $A$ to deduce information about linear degenerations of flag variety. We recall now some relevant results about quiver representations that we are going to use extensively in the next sections. 
A fundamental theorem in quiver representation theory concerns the indecomposable constituents of a generic representation $M$:
\begin{theorem}\label{CompleteDec}
Let $Q$ be a finite quiver with no loops and let $M$ be a finite dimensional representation of $Q$. The $M$ is the direct sum of uniquely determined indecomposable summands.
\end{theorem}
In particular, we need a close description of $ \Hom_Q(M,N)$ and $ \Ext^1(M,N)$ when $M,N$ are indecomposable representations of quiver $A_n$.
\begin{rmk}\label{rmk:indec} By Gabriel's Theorem, the indecomposable representations of the equioriented quiver $A_n$ are all of the form $U_{i,j}$, for $1\leq i\leq j\leq n$, where $U_{i,j}$ is the representation 
$$0\rightarrow\ldots\rightarrow 0\rightarrow \mathbb{C}\stackrel{{\rm id}}{\rightarrow} \ldots \stackrel{{\rm id}}{\rightarrow}\mathbb{C}\rightarrow 0\rightarrow\ldots\rightarrow 0,$$
supported on the vertices with indices between $i$ and $j$.
\end{rmk}
Some explicit formulae for dimension of $ \Hom_Q(U_{i,j},U_{h,k})$ and $ \Ext^1(U_{i,j},U_{h,k})$ hold:
\begin{equation}\label{hom}
\dim_\C \Hom_Q(U_{i,j},U_{h,k})=
\begin{cases}
 1 & \text{if} \;h \leq i \leq k \leq j \\
 0 & \text{otherwise.}
\end{cases}
\end{equation}
\begin{equation}\label{ext}
\dim_\C \Ext^1(U_{i,j},U_{h,k})=
\begin{cases}
 1 & \text{if} \;i+1 \leq h \leq j+1 \leq k  \\
 0 & \text{otherwise.}
\end{cases}
\end{equation}
\subsection{Coefficient Quivers}
\begin{figure}[ht]
	\centering
		\begin{tikzpicture}[baseline=(current bounding box.center), scale =.85]
			\tikzstyle{point}=[circle,thick,draw=black,fill=black,inner sep=0pt,minimum width=2pt,minimum height=2pt]
			\node[] (v00) at (-1,1) {};
            \node[] (v01) at (-1,0.333) {};
			\node[] (v02) at (-1,-0.330) {};
			\node[] (v03) at (-1,-1) {};
			\draw[fill] (v00)  circle (.05);
			\draw[fill] (v01)  circle (.05);
			\draw[fill] (v02)  circle (.05);
			\draw[fill] (v03)  circle (.05);
			\node[] (v10) at (1,1) {};
			\node[] (v11) at (1,0.333) {};
			\node[] (v12) at (1,-0.330) {};
			\node[] (v13) at (1,-1) {};
			\draw[fill] (v10)  circle (.05);
			\draw[fill] (v11)  circle (.05);
			\draw[fill] (v12)  circle (.05);
			\draw[fill] (v13)  circle (.05);
			\node[] (v20) at (3,1) {};
			\node[] (v21) at (3,0.333) {};
			\node[] (v22) at (3,-0.330) {};
			\node[] (v23) at (3,-1) {};
			\draw[fill] (v20)  circle (.05);
			\draw[fill] (v21)  circle (.05);
			\draw[fill] (v22)  circle (.05);
			\draw[fill] (v23)  circle (.05);
			
			\draw[thick ] (v00) -- (v10);
			\draw[thick  ] (v01) -- (v11);
			\draw[thick  ] (v03) -- (v13);

			\draw[thick ] (v10) -- (v20);
			\draw[thick  ] (v12) -- (v22);
			\draw[thick  ] (v13) -- (v23);
		\end{tikzpicture}
			\caption{The coefficient quiver of $M^{\bf{R}}$ for $n=3$ and ${\bf{R}}=(\{2\},\{3\})$.
			}\label{fig:H}
\end{figure}
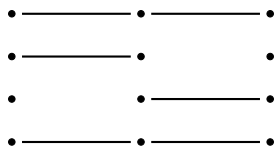

Consider now a representation $(M, F)$ of a quiver $Q$.
For every $i \in Q_0$ fix a basis $B^i=\{v_k^i\}_{k=1}^{n+1}$ of $M_i$ and set $B=\cup B_i$. 
\begin{defn}[Coefficient Quiver]\label{defn:coeffquiv}The coefficient quiver $Q(M,B)$ of $M$ with respect the basis $B$ is defined by the following data:
\begin{itemize}
\item $Q(M,B)$ has $|B|$ vertices labelled by the elements of $B$,
\item there is an arrow between $v_k^i$ and $v_h^j$ if and only if there exists an edge $\alpha \in Q_1$ between $i$ and $j$ and the coefficient of $v_h^j$ in $f_\alpha(v_k^i)$ is non zero.
\end{itemize}
\end{defn}
It is always possible to choose the basis $B$ such that to each connected component of $Q(M,B)$ corresponds an indecomposable summand of $M$ (see \cite{Ki}, Theorem 1.11).
In the following sections we will always suppose that the chosen basis $B$ used to construct the coefficient quiver $Q(M,B)$ has this property and we will omit $B$ from the notation.
\begin{example}\label{Ex:complete} Fix a basis $\{v_1, \dots, v_{n+1}\}$ of the $\C$-vector space $\C^{n+1}$. Consider the representation $M_{\bf f}$ of the quiver $A_n$ defined by the data  $M_i=\C^{n+1}$  and $f_{(i,i+1)}=\mrm{id}$ for every $i$. The coefficient quiver $Q(M_{\bf f})$ can be displayed as the union of $n+1$ linear segments of length $n$. Each linear segment corresponds to an indecomposable summand of $M$, isomorphic to $U_{1,n}$. In Figure~\ref{fig:Complete} it is displayed the described coefficient quiver for $n=3$. 
\end{example}
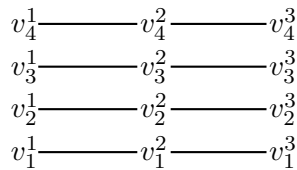
\begin{figure}[ht]
	\centering
		\begin{tikzpicture}[baseline=(current bounding box.center), scale =.85]
			\tikzstyle{point}=[circle,thick,draw=black,fill=black,inner sep=0pt,minimum width=2pt,minimum height=2pt]
			\node[] (v00) at (-0.95,1) {};
            \node[] (v01) at (-0.95,0.333) {};
			\node[] (v02) at (-0.95,-0.330) {};
			\node[] (v03) at (-0.95,-1) {};
			\node[] (v10) at (0.9,1) {};
			\node[] (v11) at (0.9,0.333) {};
			\node[] (v12) at (0.9,-0.330) {};
			\node[] (v13) at (0.9,-1) {};
			\node[] (v110) at (1.1,1) {};
			\node[] (v111) at (1.1,0.333) {};
			\node[] (v112) at (1.1,-0.330) {};
			\node[] (v113) at (1.1,-1) {};
			\node[] (v20) at (2.9,1) {};
			\node[] (v21) at (2.9,0.333) {};
			\node[] (v22) at (2.9,-0.330) {};
			\node[] (v23) at (2.9,-1) {};
			
			\node[]  at (-1,1) {$v_4^1$};
            \node[] at (-1,0.333) {$v_3^1$};
			\node[] at (-1,-0.330) {$v_2^1$};
			\node[] at (-1,-1) {$v_1^1$};
			
			\node[]  at (1,1) {$v_4^2$};
            \node[] at (1,0.333) {$v_3^2$};
			\node[] at (1,-0.330) {$v_2^2$};
			\node[] at (1,-1) {$v_1^2$};
			
			\node[]  at (3,1) {$v_4^3$};
            \node[] at (3,0.333) {$v_3^3$};
			\node[] at (3,-0.330) {$v_2^3$};
			\node[] at (3,-1) {$v_1^3$};

			\draw[thick] (v00) -- (v10);
			\draw[ thick] (v01) -- (v11);
				\draw[ thick] (v02) -- (v12);
			\draw[thick ] (v03) -- (v13);

			\draw[ thick] (v110) -- (v20);
				\draw[thick ] (v111) -- (v21);
			\draw[thick ] (v112) -- (v22);
			\draw[thick ] (v113) -- (v23);
		\end{tikzpicture}
			\caption{The coefficient quiver described in Example~\ref{Ex:complete}.
			}\label{fig:Complete}
\end{figure}
\begin{example}\label{Ex:Rdeg} Let $\bf{R}=(R_1, \dots, R_n)$ be a family of subset of $\{1, \dots, n+1\}$. In the same setting of Example~\ref{Ex:complete}, consider the representation $M^{\bf{R}}$ of the quiver $A_n$ defined by the data $M_i^{\bf{R}}=\C^{n+1}$ and $f_{(i,i+1)}=\pi_{R_i}$. Then the coefficient quiver of the representation $M^{\bf{R}}$ can be obtained from the coefficient quiver of Example~\ref{Ex:complete} deleting the edges between $v_j^i$ and $v_j^{i+1}$ if $j \in R_i$. 
In Figure~\ref{fig:H} is displayed the coefficient quiver of the representation $M^{\bf{R}}$ for ${\bf{R}}=(\{2\},\{3\})$ and $n=3$. We omitted the labeling of vertices in $Q(M^{\bf{R}})$, that must be considered as in Figure~\ref{fig:Complete}.
\end{example}
From here to the remain of the paper we deal only with representations of the equioriented quiver of type $A_n$. As mentioned before, we suppose that the coefficient quiver $Q(M)$ of a representation $M$ is always displayed with respect to a basis $B$ such that its connected components correspond to indecomposable summands of $M$. In particular, up to base change, we can suppose that there is an arrow between $v_i^j$ and $v_h^k$ \textbf{only if} $h=i$ and $k=j+1$. 
In this case we say that the (possibly disconnected)
subquiver  of $Q(M)$ spanned by the vertices associated to $\{v_i^1 \dots v_i^n \}$ is the $i$-th row of $Q(M)$.  Moreover, we will often deal with suitable subquivers of rows of $Q(M)$. We denote by $[a,b]^i_M$ the (possibly disconnected) segment in the $i$-th row of $Q(M)$ spanned by the vertices~$\{v_i^a \dots v_i^b \}$.
\begin{rmk}
We recall that, up to base change, every linear degeneration  of the flag variety ${\mathcal{F}l}(\C^{n+1})$ is $GL_{n+1}(\C)^n$-conjugate to one of the form ${\mathcal{F}l}^\mathbf{R}(\C^{n+1} )$.  Such a degeneration is isomorphic to the quiver Grassmanian $\mrm{Gr}_{\bf e}(M^{\bf{R}})$, where ${\bf e}=(1,\dots, n)$ and $M^{\bf{R}}$ is the representation described in Example~\ref{Ex:Rdeg}. As a consequence, to each $GL_{n+1}(\C)^n$-orbit $O$ it is possible to attach a coefficient quiver $Q_O$, unique up to base change. By abuse of notation, if $X$ is a degeneration in the orbit $O$, we will refer to $Q_O$ as \emph{the coefficient quiver of $X$}. Moreover, without loss of generality we can suppose that the vertices of $i$-th column of $Q_O$ are labelled from bottom to top by the vectors $v_1^i, \dots v_{n+1}^i$. We adopt this convention in everyone of the figures appearing in the remain of the paper so we omit the vertex labeling. 
\end{rmk}
\begin{example}
The coefficient quiver in Figure~\ref{fig:H} is the coefficient quiver of Feigin degeneration for~$n=3$.
\end{example}
 \begin{example} Set $\mathbf{R}_1=(\{1,2\},\{3,4\})$ and $\mathbf{R}_2=(\emptyset,\{2,3,4\})$. For a fixed basis $\{v_1, v_2, v_3, v_{4}\}$ of $\C^4$, the  degenerations ${\mathcal{F}l}^{\mathbf{R}_1}(\C^4)$ and  ${\mathcal{F}l}^{\mathbf{R}_2}(\C^4)$ are not flat because $\mathbf{R}_1$ and $\mathbf{R}_2$ do not satisfy the flatness conditions of Theorem~\ref{thm:flat}. Their coefficient quivers are displayed in Figure~\ref{fig:coeffnonflat}. 
\end{example}
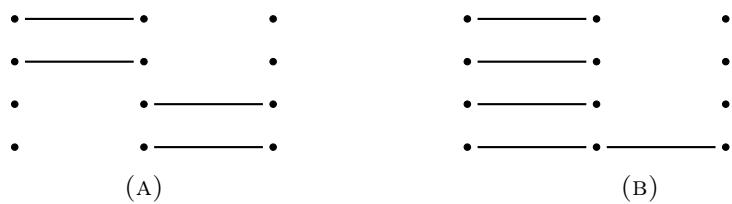
\begin{figure}[ht]
\centering
\begin{subfigure}{.3\textwidth}
\centering
		\begin{tikzpicture}[baseline=(current bounding box.center), scale =.85]
			\tikzstyle{point}=[circle,thick,draw=black,fill=black,inner sep=0pt,minimum width=2pt,minimum height=2pt]
						\node[] (v00) at (-1,1) {};
            \node[] (v01) at (-1,0.333) {};
			\node[] (v02) at (-1,-0.330) {};
			\node[] (v03) at (-1,-1) {};
			\draw[fill] (v00)  circle (.05);
			\draw[fill] (v01)  circle (.05);
			\draw[fill] (v02)  circle (.05);
			\draw[fill] (v03)  circle (.05);
			\node[] (v10) at (1,1) {};
			\node[] (v11) at (1,0.333) {};
			\node[] (v12) at (1,-0.330) {};
			\node[] (v13) at (1,-1) {};
			\draw[fill] (v10)  circle (.05);
			\draw[fill] (v11)  circle (.05);
			\draw[fill] (v12)  circle (.05);
			\draw[fill] (v13)  circle (.05);
			\node[] (v20) at (3,1) {};
			\node[] (v21) at (3,0.333) {};
			\node[] (v22) at (3,-0.330) {};
			\node[] (v23) at (3,-1) {};
			\draw[fill] (v20)  circle (.05);
			\draw[fill] (v21)  circle (.05);
			\draw[fill] (v22)  circle (.05);
			\draw[fill] (v23)  circle (.05);
			
			\draw[thick] (v00) -- (v10);
			\draw[thick] (v01) -- (v11);

			\draw[thick  ] (v12) -- (v22);
			\draw[thick ] (v13) -- (v23);
		\end{tikzpicture}
			\caption{}
			\end{subfigure}\hspace{0.1\textwidth}
			\begin{subfigure}{.3\textwidth}
		\begin{tikzpicture}[baseline=(current bounding box.center), scale =.85]
			\tikzstyle{point}=[circle,thick,draw=black,fill=black,inner sep=0pt,minimum width=2pt,minimum height=2pt]
						\node[] (v00) at (-1,1) {};
            \node[] (v01) at (-1,0.333) {};
			\node[] (v02) at (-1,-0.330) {};
			\node[] (v03) at (-1,-1) {};
			\draw[fill] (v00)  circle (.05);
			\draw[fill] (v01)  circle (.05);
			\draw[fill] (v02)  circle (.05);
			\draw[fill] (v03)  circle (.05);
			\node[] (v10) at (1,1) {};
			\node[] (v11) at (1,0.333) {};
			\node[] (v12) at (1,-0.330) {};
			\node[] (v13) at (1,-1) {};
			\draw[fill] (v10)  circle (.05);
			\draw[fill] (v11)  circle (.05);
			\draw[fill] (v12)  circle (.05);
			\draw[fill] (v13)  circle (.05);
			\node[] (v20) at (3,1) {};
			\node[] (v21) at (3,0.333) {};
			\node[] (v22) at (3,-0.330) {};
			\node[] (v23) at (3,-1) {};
			\draw[fill] (v20)  circle (.05);
			\draw[fill] (v21)  circle (.05);
			\draw[fill] (v22)  circle (.05);
			\draw[fill] (v23)  circle (.05);
			
			\draw[thick ] (v00) -- (v10);
			\draw[thick] (v01) -- (v11);
			\draw[thick] (v02) -- (v12);
			\draw[thick  ] (v03) -- (v13);

			\draw[thick  ] (v13) -- (v23);
		\end{tikzpicture}
			\caption{} 
			\end{subfigure}
			\caption{Coefficient quivers of two non flat degenerations of ${\mathcal{F}l}(\C^{4})$}\label{fig:coeffnonflat}
\end{figure}
\section{Torus Actions and Cohomology}\label{sec:3} In recent papers \cite{MA} and \cite{MAII},  it is proved that every linear degeneration of the flag variety can be
endowed with a structure of GKM variety, under the action of a suitable algebraic torus~$T$. Their results hold in the more general context of cyclic quiver Grassmanians and as a consequence of this more general point of view, it is possible to achieve a description of the moment graph in a purely combinatorial way using quiver representation theory. We used these  techniques to prove our smoothness criteria. In this section we recall some basic facts about GKM varieties and the main results contained in \cite{MA}
 and \cite{MAII}. \subsection{GKM varieties and their cohomology.}
Let $X$ be a complex projective algebraic $T$-variety, i.e. $X$ is equipped by an action of an algebraic torus $T \simeq (\C^*)^r$ (for some integer $r$). We denote by $H^*_T(X)$ the $T$ equivariant cohomology ring of $X$ (with coefficients in $\Q$). 
\begin{defn}\label{defn:GKM} A $T$ variety $X$ is a $GKM$ variety if: 
\begin{enumerate}
    \item The number of $T$-fixed points and of $1$ dimensional $T$-orbits in $X$ is finite
    \item The usual cohomology of $X$ can be recovered by ${H_T^*}(X)$ by scalar extension: 
    \[H^*(X) \simeq H^*_T(X) \otimes_{H^*_T(pt)} \Q.\]
\end{enumerate}
\end{defn}
\begin{rmk}In our context the latter condition is equivalent to require that rational cohomology of $X$ vanishes in odd degrees. (c.f.r. \cite{MA}, Section 1)
\end{rmk} 
\begin{defn}
Let $X$ be a $T$ variety and $\chi$ a cocharacter of $T$. We say that $\chi$ is a \emph{generic cocharacter (for the $T$ action on $X$)} if $X^{\chi(\C^*)}=X^T$.
\end{defn}
Existence of a generic cocharacter implies the existence of a nice cell decomposition of $X$.
\begin{defn}[Białynicki-Birula decomposition]
Let $X$ be an algebraic variety equipped with an action of the one dimensional algebraic torus $\C^*$. Denote by $X_1, \dots, X_n$ the components of $X^T$ and set \[C_i=\{x \in X \; | \; \lim_{z \rightarrow 0} z \cdot x \in X_i\}.\]
Then $X$ is the disjoint union of $C_1, \dots, C_n$.
\end{defn}
Consider $x\in X^T$, we will denote by $C(x)$ the unique attractive cell $C_i$ in the Białynicki-Birula decomposition of $X$ that contains $x$.
\begin{defn}[Moment Graph]\label{defn:moment}
Let $X$ be a $T$ variety. Its Moment Graph is the graph $G(X,T)$ defined by the following data:
\begin{itemize}
\item The vertices of $G(X,T)$ are indexed by the set of fixed points $X^T$,
\item There is an (oriented) edge from $x$ to $y$ if $x \in \overline{C(y)}$
\end{itemize}
\end{defn}
One of the main results in Goresky - Kottwitz - MacPherson theory (c.f.r. \cite{GKM}) assert that the $T$-equivariant cohomology of a a GKM variety can be explicitly described using its moment graph. 

We say that a cell $C_i$ defined as above is \emph{rational} if $C_i$ is rationally smooth at every $x \in C_i$, i.e. for every $x \in C_i$, the top rational cohomology group  $H^{2 \dim_\C C_i}(C_i,C_i \setminus \{x\})$ has rank 1 and  $H^{i}(C_i,C_i \setminus \{x\})=0$ for every $i\neq 2 \dim_\C C_i$. 
\begin{defn}[BB filterable variety]
A $T$ variety is BB filterable if:
\begin{enumerate}
    \item[BB1)] the set $X^T$ is finite,
    \item[BB2)] $X$ admit a BB decomposition in rational cells, induced by the action a suitable generic cocharacter $\chi$ of $T$.
\end{enumerate}
\end{defn}
In \cite{MA} it is proved that, under the action of a suitable algebraic torus $T$,  cyclic quiver Grassmanians have a structure of BB-filterable varieties.
As a consequence of (\cite{MA}, Theorem 1.14) a projective BB filterable $T$-variety is GKM and consequently linear degenerations of the flag variety of type $A$ are GKM. 
\subsection{Tori Actions on Linear Degenerations.}
\begin{figure}[ht]
\centering
\begin{subfigure}{.3\textwidth}
\centering
		\begin{tikzpicture}[baseline=(current bounding box.center), scale =.7]
			\tikzstyle{point}=[circle,thick,draw=black,fill=black,inner sep=0pt,minimum width=2pt,minimum height=2pt]
			\node[] (v00) at (-1,1) {};
            \node[] (v01) at (-1,0.333) {};
			\node[] (v02) at (-1,-0.330) {};
			\node[] (v03) at (-1,-1) {};
			\draw[fill,red] (v00)  circle (.07);
			\draw[fill] (v01)  circle (.05);
			\draw[fill] (v02)  circle (.05);
			\draw[fill] (v03)  circle (.05);
			\node[] (v10) at (1,1) {};
			\node[] (v11) at (1,0.333) {};
			\node[] (v12) at (1,-0.330) {};
			\node[] (v13) at (1,-1) {};
			\draw[fill,red] (v10)  circle (.07);
			\draw[fill,red] (v11)  circle (.07);
			\draw[fill] (v12)  circle (.05);
			\draw[fill] (v13)  circle (.05);
			\node[] (v20) at (3,1) {};
			\node[] (v21) at (3,0.333) {};
			\node[] (v22) at (3,-0.330) {};
			\node[] (v23) at (3,-1) {};
			\draw[fill,red] (v20)  circle (.07);
			\draw[fill] (v21)  circle (.05);
			\draw[fill,red] (v22)  circle (.07);
			\draw[fill,red] (v23)  circle (.07);

			\draw[ultra thick, red] (v00) -- (v10);
			\draw[,  ] (v01) -- (v11);
			\draw[,  ] (v03) -- (v13);

			\draw[ultra thick, red] (v10) -- (v20);
			\draw[,  ] (v12) -- (v22);
			\draw[,  ] (v13) -- (v23);
		\end{tikzpicture}
			\caption{$S=(\{4\}, \{4,3\}, \{4,2,1\})$}\label{fig:SCSa}
\end{subfigure}\hspace{0.2\textwidth}
\begin{subfigure}{.3\textwidth}
\centering
		\begin{tikzpicture}[baseline=(current bounding box.center), scale =.7]
			\tikzstyle{point}=[circle,thick,draw=black,fill=black,inner sep=0pt,minimum width=2pt,minimum height=2pt]
			\node[] (v00) at (-1,1) {};
            \node[] (v01) at (-1,0.333) {};
			\node[] (v02) at (-1,-0.330) {};
			\node[] (v03) at (-1,-1) {};
			\draw[fill] (v00)  circle (.05);
			\draw[fill] (v01)  circle (.05);
			\draw[fill] (v02)  circle (.05);
			\draw[fill,red] (v03)  circle (.07);
			\node[] (v10) at (1,1) {};
			\node[] (v11) at (1,0.333) {};
			\node[] (v12) at (1,-0.330) {};
			\node[] (v13) at (1,-1) {};
			\draw[fill,red] (v10)  circle (.07);
			\draw[fill] (v11)  circle (.05);
			\draw[fill] (v12)  circle (.05);
			\draw[fill, red] (v13)  circle (.07);
			\node[] (v20) at (3,1) {};
			\node[] (v21) at (3,0.333) {};
			\node[] (v22) at (3,-0.330) {};
			\node[] (v23) at (3,-1) {};
			\draw[fill,red] (v20)  circle (.07);
			\draw[fill,red] (v21)  circle (.07);
			\draw[fill,red] (v22)  circle (.07);
			\draw[fill] (v23)  circle (.05);

			\draw[, ] (v00) -- (v10);
			\draw[ultra thick, red ] (v03) -- (v13);
			\draw[,  ] (v01) -- (v11);

			\draw[ultra thick, red] (v10) -- (v20);
			\draw[,  ] (v12) -- (v22);
			\draw[,  ] (v13) -- (v23);
		\end{tikzpicture}
		\caption{$S'=(\{1\}, \{1,4\}, \{4,3,2\})$}\label{fig:SCSb}
\end{subfigure}
\caption{Two subquivers and their associated families of subsets}\label{fig:SCS}
\end{figure}
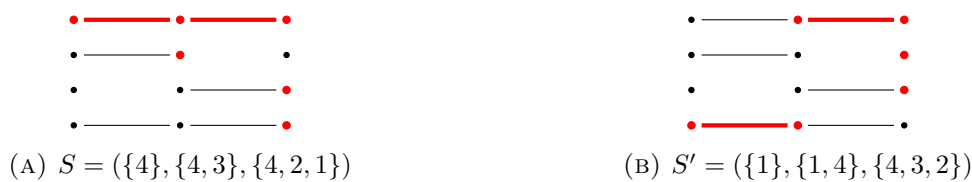

In this section we recall the construction of the torus action studied in \cite{MA} and \cite{MAII}. This action is inspired by the one used in \cite{CI} to compute the Euler characteristic of certain quiver Grassmaniains. Moreover, in \cite{MA} it is proved that a generic chocaracter $\chi$ for this action exists and the associated $\C^*$-action induces a Bialynicki-Birula decomposition. \\
Let $M^{\bf{R}}$ be the representation of the quiver $A_n$ associated to ${\mathcal{F}l}^\mathbf{R}(\C^{n+1})$ as described in Example~\ref{Ex:Rdeg}. Each indecomposable summand of $M^{\bf{R}}$ is of the form $U_{h,k}$ (Remark~\ref{rmk:indec}) and can be identified with a maximal connected component of some row $r_i$ of $Q(M^{\bf{R}})$. In our notation we stress the fact that an indecomposable $U$ isomorphic to $U_{h.k}$ is identified with a connected component of $r_i$ denoting it by $U_{h,k}^i$. 
Suppose now that $M^{\bf{R}}$ has $d_0$ indecomposable summands, order them and consider the torus $T=(\C^*)^{d_0}$. If $U_{h,k}^i$ is the $j$-th indecomposable summand, with respect to the chosen ordering, we define the $T$ action on $U_{h,k}^i$ by the following formula: let $\gamma=(\gamma_1, \dots, \gamma_{d_0})$ be the generic element of $T$, then for every $t \in \{h, \dots, k\}$ set
\begin{equation}\label{eqn:action}\gamma \cdot v^i_t = \gamma_j v^i_t. \end{equation}
We want now to define a generic cocharacter for this action.
\begin{defn}[Attractive Grading, \cite{MA}, Definition 3.10] An attractive grading on $Q(M)_0$ is a map $\mrm{wt}: Q(M)_0 \rightarrow \Z$ that satisfies the following properties: 
\begin{enumerate}
    \item For every $i \in Q_0$ it holds that $\mrm{wt}(v_k^i)>\mrm{wt}(v_h^i)$, whenever $k>h$;
    \item For every edge $\alpha \in Q_1$ there exists $d(\alpha) \in \Z$ such that $\mrm{wt}(v_k^{t(\alpha)})>\mrm{wt}(v_h^{s(\alpha)}) + d(\alpha)$ whenever the is and ordinted edge $\mrm{wt}(v_h^{s(\alpha)})\rightarrow \mrm{wt}(v_k^{t(\alpha)})$ in $Q(M)$.
\end{enumerate}
\end{defn}
Fix an attractive grading on the vertices of $Q(M^{\bf{R}})$.
If $U^i_{h,k}$ is the $j$-th indecomposable summand of $M^{\bf{R}}$, with respect to the ordering of indecomposable that we have chosen to define the $T$ action, let us denote by $b_j$ the starting vertex of the connected component associated to $U^i_{h,k}$ in $Q(M^{\bf{R}})$. The  map $\chi : \C^* \rightarrow T$ defined by the assignment 
\begin{equation}\label{eqn:cochar} z \in \C^* \rightarrow (z^{\mrm{wt}(b_1)}, \dots, z^{\mrm{wt}(b_{d_0})} )\end{equation}
defines a cocharacted of $T$. The following Theorems are proved in \cite{MA}:
\begin{thm}[\cite{MA}, Theorem 5.14]
Let $\mrm{wt}$ be an attractive grading, $T= (\C)^{d_0}$ the torus acting on $\mrm{Gr}_{\bf e}(M^{\bf{R}})$ by Formula~\ref{eqn:action} and $\chi$ the cocharacter defined by (\ref{eqn:cochar}). Then $\chi$ is a generic cocharacter for the action of $T$ on $\mrm{Gr}_{\bf e}(M^{\bf{R}})$
\end{thm}
\begin{thm}[\cite{MA}, Proposition 6.4]
The quiver Grassmanian $\mrm{Gr}_{\bf e}(M^{\bf{R}})$ is a $BB$-filterable projective $T$-variety. Moreover, the $T$-action induces on $\mrm{Gr}_{\bf e}(M^{\bf{R}})$ a structure of $GKM$ variety. 
\end{thm}
Consequently, some natural questions arise, in particular determining a description of
the moment graph  for the $T$ action on ${\mathcal{F}l}^\mathbf{R}(\C^{n+1})$ in a combinatorial flavour results to be very interesting. In (\cite{MA}, Section 6.2) it is proved that fixed points for the $T$ action can be identified with suitable subquivers of the coefficient quiver $Q(M^{\bf{R}})$.
 \begin{defn}[Successor Closed Subquiver, c.f.r.\cite{MA}, Definition 6.7]\label{defn:scs}  A successor closed subquiver of $Q(M^{\bf{R}})$ is a full subquiver $Q'$ of $Q(M^{\bf{R}})$ such that 
 \begin{itemize}
     \item if $v \in Q'_0$ and $v=s(\alpha)$ for some $\alpha \in Q(M^{\bf{R}})_1$, then $t(\alpha) \in Q'_0$
     \item $|Q'_0 \cap B_i|=i$ for all $i \leq n$.
     \end{itemize}
 \end{defn}
  \begin{theorem}[c.f.r. \cite{MA},Theorem 6.15]
 The points in $X^T$ are in bijection with successor closed subquivers of $ Q(M^{\bf{R}})$.
 \end{theorem}
Successor closed subquivers of $Q(M^{\bf{R}})$ can be described using certain sequences of integers. 
 \begin{defn}[${\bf R}$-Admissible Sequence]\label{defn:admissible}
 A family  $S=(S_1, \dots, S_{n})$ of subsets $S_i \subset \{1 \dots, n+1\}$ is $\mathbf{R}$-admissible sequence if $|S_i|=i$ for every $i \leq n$ and $S_{i} \subset S_{i+1}\cup R_i $.
 \end{defn}
The condition $S_{i} \subset S_{i+1}\cup R_i $ implies that any $\mathbf{R}$-admissible sequence $S$ define a successor closed subquiver of $Q(M^{\bf{R}})$, and \emph{vice versa} every successor closed subquiver of $Q(M^{\bf{R}})$ is associated to a unique $\mathbf{R}$-admissible sequence. We denote  the subquiver associated to $S$ by $Q_S$ and by $p_S$ the corresponding point in $X^T$.
\begin{example}
We recall that for Feigin Degeneration  for every $i$ we have $R_i=\{i+1\}$. A fixed point $p_S$ is then associated to an $\mathbf{R}$-admissible family $S=(S_1, \dots, S_n)$ such that $S_i \subset S_{i+1} \cup \{i+1\}$. In Figure~\ref{fig:SCS} are displayed two subquivers of the coefficient quiver of the Feigin degeneration for $n=3$. Both of them are encoded by families of subsets of $\{1,2,3,4\}$. The family $S$ associated to the subquiver in Figure~\ref{fig:SCSa} satisfies the condition of being an $\mathbf{R}$-admissible sequence and, in fact, the associated subquiver is successor closed. Conversely the subquiver in Figure~\ref{fig:SCSb} is not successor closed and the sequence $S'$ is not $\mathbf{R}$-admissible.
\end{example}

\begin{figure}[ht]
\centering
\begin{subfigure}{.3\textwidth}
\centering
		\begin{tikzpicture}[baseline=(current bounding box.center), scale =.7]
			\tikzstyle{point}=[circle,thick,draw=black,fill=black,inner sep=0pt,minimum width=2pt,minimum height=2pt]
			\node[] (v00) at (-1,1) {};
            \node[] (v01) at (-1,0.333) {};
			\node[] (v02) at (-1,-0.330) {};
			\node[] (v03) at (-1,-1) {};
			\draw[fill,red] (v00)  circle (.07);
			\draw[fill] (v01)  circle (.05);
			\draw[fill] (v02)  circle (.05);
			\draw[fill] (v03)  circle (.05);
			\node[] (v10) at (1,1) {};
			\node[] (v11) at (1,0.333) {};
			\node[] (v12) at (1,-0.330) {};
			\node[] (v13) at (1,-1) {};
			\draw[fill,red] (v10)  circle (.07);
			\draw[fill,red] (v11)  circle (.07);
			\draw[fill] (v12)  circle (.05);
			\draw[fill] (v13)  circle (.05);
			\node[] (v20) at (3,1) {};
			\node[] (v21) at (3,0.333) {};
			\node[] (v22) at (3,-0.330) {};
			\node[] (v23) at (3,-1) {};
			\draw[fill,red] (v20)  circle (.07);
			\draw[fill,red] (v21)  circle (.07);
			\draw[fill,red] (v22)  circle (.07);
			\draw[fill] (v23)  circle (.05);

			\draw[ultra thick, red] (v00) -- (v10);
			\draw[,  ] (v01) -- (v11);
			\draw[,  ] (v03) -- (v13);

			\draw[ultra thick, red] (v10) -- (v20);
			\draw[,  ] (v12) -- (v22);
			\draw[,  ] (v13) -- (v23);
		\end{tikzpicture}
		\caption{$S=(\{4\}, \{3,4\}, \{2,3,4\})$}\label{mutA}
\end{subfigure}\hspace{0.2\textwidth}
\begin{subfigure}{.3\textwidth}
\centering
		\begin{tikzpicture}[baseline=(current bounding box.center), scale =.7]
			\tikzstyle{point}=[circle,thick,draw=black,fill=black,inner sep=0pt,minimum width=2pt,minimum height=2pt]
			\node[] (v00) at (-1,1) {};
            \node[] (v01) at (-1,0.333) {};
			\node[] (v02) at (-1,-0.330) {};
			\node[] (v03) at (-1,-1) {};
			\draw[fill] (v00)  circle (.05);
			\draw[fill,red] (v01)  circle (.07);
			\draw[fill] (v02)  circle (.05);
			\draw[fill] (v03)  circle (.05);
			\node[] (v10) at (1,1) {};
			\node[] (v11) at (1,0.333) {};
			\node[] (v12) at (1,-0.330) {};
			\node[] (v13) at (1,-1) {};
			\draw[fill,red] (v10)  circle (.07);
			\draw[fill,red] (v11)  circle (.07);
			\draw[fill] (v12)  circle (.05);
			\draw[fill] (v13)  circle (.05);
			\node[] (v20) at (3,1) {};
			\node[] (v21) at (3,0.333) {};
			\node[] (v22) at (3,-0.330) {};
			\node[] (v23) at (3,-1) {};
			\draw[fill,red] (v20)  circle (.07);
			\draw[fill,red] (v21)  circle (.07);
			\draw[fill,red] (v22)  circle (.07);
			\draw[fill] (v23)  circle (.05);

			\draw[, ] (v00) -- (v10);
			\draw[ultra thick, red ] (v01) -- (v11);
			\draw[,  ] (v03) -- (v13);

			\draw[ultra thick, red] (v10) -- (v20);
			\draw[,  ] (v12) -- (v22);
			\draw[,  ] (v13) -- (v23);
		\end{tikzpicture}
		\caption{$S'=(\{3\}, \{3,4\}, \{2,3,4\})$}\label{MutB}
\end{subfigure}
\caption{A Mutation from the $4$-th row to the $3$-rd}\label{mutation}
\end{figure}
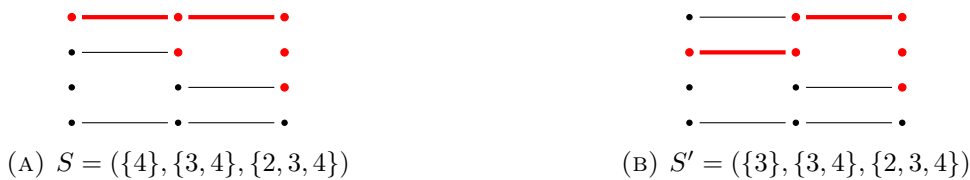

Aiming to clarify how combinatorics of admissible sequences plays a very relevant role in studying the singularity of fixed points, we think it could be very useful now to remark that for each $\mathbf{R}$-admissible sequence $S$ we have three distinct objects: a fixed point for the torus action $p_S$, a subrepresentation $M_S$ of $M^{\bf{R}}$ and a successor closed subquiver $Q_S$ of $Q(M^{\bf{R}})$. We explain now how these three objects are linked.\\
A crucial fact to prove our smoothness criteria links $M_S$ to the dimension of tangent space to~${\mathcal{F}l}^\mathbf{R}(\C^{n+1})$ at~$p_S$:
\begin{theorem}[c.f.r. \cite{CR}, Proposition 6]\label{HomTg}
\[\dim \Hom(M_S, M^{\bf{R}}/M_S) = \dim T_{p_S}({\mathcal{F}l}^\mathbf{R}(\C^{n+1}))\]
\end{theorem}
\begin{rmk}\label{ExtZero}
The indecomposable summands of quotient representation $M^{\bf{R}}/M_S$ corresponds to the connected component of the Coefficient Quiver $Q(M^{\bf{R}})\setminus Q_S$. In particular if $M_S$ has dimension vector $(1, \dots, n)$, then $M^{\bf{R}}/M_S$ has dimension vector $(n, \dots, 1)$ and then $$\langle \dim M_S , \dim M^{\bf{R}}/M_S \rangle = \frac{n(n+1)}{2}.$$ This implies, as a consequence of Formula~(\ref{ExtHom}) and of Theorem \ref{HomTg}, that if ${\mathcal{F}l}^\mathbf{R}(\C^{n+1})$ is flat, a point $p_S$ is smooth if and only if $\Ext^1(M_S, M^{\bf{R}}/M_S)=0$.
\end{rmk}
Let us denote by $SC(M^{\bf{R}})$ the set of successor closed subquivers of $Q(M^{\bf{R}})$.
It is proved in \cite{MA} that an attractive grading on the vertices of the  coefficient quiver $Q(M^{\bf{R}})$ induces an ordering on the elements of $SC(M^{\bf{R}})$ that is isomorphic to the partial order on fixed points in the moment graph of ${\mathcal{F}l}^{\bf{R}}(\C^{n+1})$. It is possible to give a combinatorial description of this ordering in terms of \emph{mutations} of a successor closed subquiver $Q_S$.\\
We remark that, by our choice of the basis $B$, each connected component of $Q(M^{\bf{R}})$ is a linear segment, i.e. the valence of each vertex is smaller or equal than 2.
\begin{defn}[Movable part of $Q_S$, c.f.r. \cite{MA}, Definition 6.8]
A movable part of a linear segment $L \subset Q_S $ is a connected subquiver $L' \subset L$ such that $L'$ has the same starting vertex of $L$.
\end{defn}
In other words, a movable part is a connected segment $[a,b]_{Q_S}^i$ of $Q_S$ such that $[a-1,b]_{Q_S}^i$ is not a connected linear segment of $Q_S$.
In (\cite{MA}, Section 6) the authors describe combinatorially the covering relations in the moment graph in the following way: let $Q_S$ and $Q_{S'}$ be successor closed subquivers of $Q(M^{\bf{R}})$, we have that $Q_S $ is covered by $ Q_{S'}$ if and only if there exists a \emph{fundamental mutation}(\cite{MA}, Definition 6.8) from $Q_{S'}$ to $Q_{S}$, i.e. it is possible to obtain $Q_{S}$ from $Q_{S'}$ moving down \emph{exactly one} movable part of $Q_{S'}$ with respect to a fixed attractive ordering on the elements of $B$.\\
One of the aims of this note is to link mutations with the property of a fixed point in a flat degeneration of being singular.
In particular, we are interested only to the valence of vertices in the moment graph, so we relax the definition of mutations as presented in \cite{MA}:
\begin{defn}[Mutation]\label{defn:mutation}
Consider $Q_S, Q_{S'} \in SC(M^{\bf{R}})$, we say that there is a mutation from $Q_{S'} $ to $Q_S$ if $Q_S$ is obtained from $Q_{S'}$ moving \emph{exactly} one movable part.
\end{defn}
In other words, if $Q_S, Q_{S'} $ are successor closed subquivers of $Q(M^{\bf{R}})$, there is a mutation from $Q_{S'} $ to $Q_S$ if and only if $Q_S \setminus Q_{S'}$ and $Q_{S'} \setminus Q_{S}$ are both movable parts of the same length (i.e. are spanned by the same number of vertices).
In terms of connected sub-segments of the coefficient quiver $Q(M^{\bf{R}})$ we say that $Q_S$ is obtained from $Q_{S'}$ by a mutation if there exists two connected segment $I=[a,b]^i_{Q_{S'}}$ and $J=[a,b]^j_{Q(M^{\bf{R}})}$ such that $Q_S=(Q_{S'} \cup J) \setminus I $.
We denote by $\mrm{Mut}(Q_S)$ the set of mutations of a successor closed subquiver $Q_S\in SC(M^{\bf{R}})$.
\begin{rmk}\label{rmk:valencemutations}By definition of mutation and as a consequence of (\cite{MA}, Theorem 6.15), the cardinality of $\mrm{Mut}(Q_S)$ equals the valence of $p_S$ in the moment graph.\end{rmk}
\begin{rmk}\label{segmentmutations}We have a mutation from $I=[a,b]^i_{Q_S}$ to $J=[a,b]^j_{Q(M^{\bf{R}})}$ only if $j \in S_{b+1}$ or if $f_b(j)=0$ in $Q(M^{\bf{R}})$. In this case we write $[a,b]^i_{Q_S}\rightarrow j$.\end{rmk}
\begin{example} In Figure~\ref{mutA} is displayed a successor closed subquiver $Q_S$ of the coefficient quiver of Feigin Degeneration for $n=3$. The successor closed subquiver $Q_{S'}$ in Figure~\ref{MutB} is obtained from the previous one by the mutation $[1]^4_{Q_S} \rightarrow 3$.\end{example}

We recall that successor closed subquivers of $Q(M^{\bf{R}})$ are in bijections with $T$-fixed point in  ${\mathcal{F}l}^\mathbf{R}(\C^{n+1})$ and that each of these points is indexed by an admissible sequence $S=(S_1, \dots, S_n)$ of subsets of $[\{1, \dots, n+1\}$. By abuse of notation, in this case we will denote by $[a,b]^i_S$ the generic movable part in $Q_S$ and by $\mrm{Mut}(p_S)$ the set of possible mutations of the successor closed subquiver associated to the point $p_S$. 
\begin{proposition}\label{ValSmooth} Let $p_S$ be a fixed point, then $|\mrm{Mut}(p_S)|= \mrm{dim}T_{p_S}{\mathcal{F}l}^\mathbf{R}(\C^{n+1})$. \end{proposition} 
\proof Because of Theorem~\ref{HomTg} and by Formula~(\ref{hom}), it is enough to prove that there is a mutation $[a,b]^i_{S}\rightarrow j$ if and only if there exist and indecomposable summand $U_{a,b}^i$ of $M_S$ and an indecomposable summand $U_{a',b}^j$ of $M^{\bf{R}} / M_S$ such that  $\Hom (U_{a,b}^i, U_{a',b}^j) \neq 0$. In is clear by definition that if we have a mutation $[a,b]^i_{S}\rightarrow j$, then $U_{a,b}$ is an indecomposable summand of $M_S$ (corresponding to the connected segment $[a,b]^i_{S}$ of $Q_S$). Moreover, by definition, the segment $[a,b]^j_{Q(M^{\bf{R}})}$ is is connected. By Remark~\ref{segmentmutations}, a maximal connected segment of the $j$-th row of  $Q(M^{\bf{R}}) \setminus Q_S$ is of the form $[a',b]^j$ with $a' \leq a$. As a consequence, the linear segment  $[a',b]^j$ corresponds to an indecomposable summand $U_{a',b}$ of $M^{\bf{R}} / M_S$ and $\Hom (U_{a,b}, U_{a',b}) \neq 0$. Observe that such a pair of indecomposables is uniquely determined by the mutation $[a,b]^i_{S}\rightarrow j$. On the other hand, if $\Hom (U_{a,b}^i, U_{a',b'}^j) \neq 0$, where $U_{a,b}^i$ is an indecomposable summand of $M_S$ and  $U_{a',b'}^j$ is an indecomposable summand of $M^{\bf{R}} / M_S$, then by Formula~(\ref{hom}) we have $a' \leq a \leq b' \leq b$. Consequently $[a,b']^i_S$ is a a connected segment of the $i$-th row of $Q_S$ such that $[a-1,b']^i_S$ is not a linear segment of $Q_S$. Moreover $[a',b']^i_S$ is a maximal connected segment of $Q(M^{\bf{R}}) \setminus Q_S$ and then the pair $(U_{a,b}^i, U_{a',b'}^j)$ corresponds to an uniquely determined mutation $[a,b']^i_{S}\rightarrow j$.
\endproof
In particular, observe that if $|\mrm{Mut}(p_S)|=\binom{n+1}{2}$, then $p_S$ is \emph{always} smooth. As a consequence of Remark~\ref{rmk:valencemutations}, it is possible to link the dimension of the tangent plane at $p_S$ to the valence of $p_S$ in the moment graph.
\begin{corollary}\label{SingDim} The valence of the vertex corresponding to $p_S$ in the moment graph is equal to the dimension of tangent space to ${\mathcal{F}l}^\mathbf{R}(\C^{n+1})$ at $p_S$. 
\end{corollary}

\section{A Combinatorial Formula for the Dimension of Tangent Space}\label{sec:formula}
The first step for our combinatorial characterization of smooth fixed points is to define a suitable set that encodes relevant informations about the dimension of the tangent space at a fixed point $p_S$.
\begin{defn}
Consider the family ${\bf R}=(R_1, \dots, R_{n-1})$ of subsets of $\{1, \dots, n+1\}$. We define the set of $i$-positions of $\bf{R}$ as the set 
\[\pos (i) = \{j \in \{1, \dots, n-1\}\; | \; i \in R_j\}.\] 
\end{defn}
We will also refer to $\pos (i) $ as the set of positions of the $i$-th row in~$Q(M^{\bf R})$.
\begin{defn}\label{SingSet} We define $\mathrm{Sing}_i (p_S)$ as the set of pairs $(j,h+1)$, with $h \in \pos(i)$ and $j \notin S_{h+1}$,  such that there exists $ k \leq h, \, k \in \pos(j)$ with the following characteristics:
\begin{enumerate}
    \item  $i \in S_t$ for all $ t $  such that $k \leq t \leq h$
    \item the segment spanned by the vertices $\{v_i^k, \dots , v_i^h\}$ is connected in $Q_S$.
 \end{enumerate}
\end{defn}
 We say that an element of $\mathrm{Sing}_i (p_S)$ is a singular point for the $i$-th row of $Q_S$.
 \begin{defn}The set of singularities of $Q_S$ is the set $\mathrm{Sing} (p_S)= \bigsqcup \mathrm{Sing}_i (p_S)$.
 \end{defn}
 \begin{example}\label{ex:nonsmoothF3}
We remark that a pair $(h,k) \in \mathrm{Sing} (p_S)$ can be displayed as a vertex of $Q(M^{\bf R})$ (with some multiplicities). In Figure~\ref{fig:sing3} the set $\mathrm{Sing}(p_S)$ has a unique contribution coming from singular points for $3$-rd row. The set $\mathrm{Sing}_3(p_S)$ is composed by the pair of indices $(3,2)$ that corresponds to the blue vertex in Figure~\ref{fig:sing3}.\end{example}
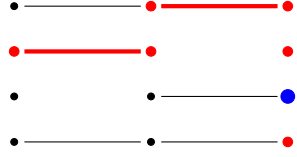
\begin{figure}[ht]
\centering
		\begin{tikzpicture}[baseline=(current bounding box.center), scale =.9]
			\tikzstyle{point}=[circle,thick,draw=black,fill=black,inner sep=0pt,minimum width=2pt,minimum height=2pt]
			\node[] (v00) at (-1,1) {};
            \node[] (v01) at (-1,0.333) {};
			\node[] (v02) at (-1,-0.330) {};
			\node[] (v03) at (-1,-1) {};
			\draw[fill] (v00)  circle (.05);
			\draw[fill,red] (v01)  circle (.07);
			\draw[fill] (v02)  circle (.05);
			\draw[fill] (v03)  circle (.05);
			\node[] (v10) at (1,1) {};
			\node[] (v11) at (1,0.333) {};
			\node[] (v12) at (1,-0.330) {};
			\node[] (v13) at (1,-1) {};
			\draw[fill,red] (v10)  circle (.07);
			\draw[fill,red] (v11)  circle (.07);
			\draw[fill] (v12)  circle (.05);
			\draw[fill] (v13)  circle (.05);
			\node[] (v20) at (3,1) {};
			\node[] (v21) at (3,0.333) {};
			\node[] (v22) at (3,-0.330) {};
			\node[] (v23) at (3,-1) {};
			\draw[fill,red] (v20)  circle (.07);
			\draw[fill,red] (v21)  circle (.07);
			\draw[fill,blue] (v22)  circle (.1);
			\draw[fill,red] (v23)  circle (.07);

			\draw[] (v00) -- (v10);
			\draw[ultra thick, red] (v01) -- (v11);
			\draw[,  ] (v03) -- (v13);

			\draw[ultra thick, red] (v10) -- (v20);
			\draw[,  ] (v12) -- (v22);
			\draw[,  ] (v13) -- (v23);
		\end{tikzpicture}
		\caption{The set $\mathrm{Sing}_3(p_S)$, where~$S=\{\{3\},\{3,4\},\{1,3,4\}\}$}\label{fig:sing3}
\end{figure}
We prove now that, for \emph{every} linear degeneration of the form in ${\mathcal{F}l}^{\bf{R}}(\C^{n+1})$, the following formula for the dimension of the tangent space at $p_S$ holds:
\begin{equation}\label{dim}
    \dim T_{p_S}{\mathcal{F}l}^{\bf R}(\C^{n+1}) = \frac{n(n+1)}{2}+ |\mathrm{Sing} (p_S)|.
\end{equation}
More precisely, we prove that $\mathrm{Sing} (p_S)$ measures the dimension of $\Ext^1(M_S, M^{\bf R}/M_S)$.
\begin{proposition}
$$|\mathrm{Sing} (p_S)|=\dim \Ext^1(M_S, M^{\bf R}/M_S)$$
\end{proposition}
\proof First of all observe that $(j,h+1)\in \mathrm{Sing}_i (p_S)$ implies that there exist some uniquely determined indices $a= \max\{t\in \pos(j)| t\leq h\}+1$ and $b>h$ such that $[a,b]^j_M$  is a connected segment in $Q(M^{\bf R})\setminus Q_S$. This segment is associated to an indecomposable subrepresentation $U^j_{a,b}$ of $M^{\bf R}/M_S$.
Consider the two integers $\bar{k}$ and $k_i$ such that 
\[\bar{k}=\max \{t \in \pos(i) \cup \{1\}| t \leq h\} \qquad k_i=\min\{t \leq h | i \in S_{s} \; \forall s,  t \leq s \leq h\}.\]
The segment $[\max\{\bar{k},k_i\}, h]^i_S$ is then a connected component of the $i$-th row of $Q_S$ that corresponds to an indecomposable summand $U^i_{\max\{\bar{k},k_i\}, h}$ of $M_S$.
Moreover, the condition $2)$ in Definition~\ref{SingSet}, implies that $\max\{\bar{k},k_i\} \leq a \leq h < b$.
Consequently, for each element in $(j,h+1)\in \mathrm{Sing}_i (p_S)$ we constructed pair $(U^i_{\max\{\bar{k},k_i\}, h} , U^j_{a,b})$ of indecomposable representations in $M_S$ and $M/M_S$ respectively, such that $\Ext^1(U^i_{\max\{\bar{k},k_i\}, h} , U^j_{a,b}) \neq 0$. Observe now that, once we fixed the index $i$, the pair $(U^i_{\max\{\bar{k},k_i\}, h} , U^j_{a,b})$ is uniquely determined by $(j,h+1) \in \mrm{Sing}_i(p_S)$.
Then we need only to prove that each pair $(U_{k,h}, U_{a,b})$ of indecomposable summands of $M_S$ and $M^{\bf R}/M_S$ respectively, such that $\Ext^1(U_{k,h}, U_{a,b})\neq 0$, is associated to a pair $(j,h+1) \in \mathrm{Sing}(p_S)$.
By condition expressed in Formula (\ref{ext}), we have $k+1 \leq a \leq h+1 \leq b$. Without loss of generality we can suppose that $U_{k,h}$ and $U_{a,b}$ corresponds to connected components of the $i$-th row of $Q_S$ and of $j$-th row of $Q(M^{\bf R})\setminus Q_S$ respectively. This implies that the pair of indices $(j,h+1)$ satisfies the conditions in Definition~\ref{SingSet} and then $(j,h+1)\in \mathrm{Sing}_i (p_S)$.
\endproof
Formula (\ref{dim}) now follows from Formula (\ref{ExtHom}).
\begin{example} We provide now an example of singularity sets for two different rows of the same $Q_S$  that are not disjoint. Consider the flat degeneration ${\mathcal{F}l}^{\bf R}(\C^7)$ defined by the sequence of subsets  ${ \bf R}=\{\{5\}, \{6\}, \emptyset, \{3,6\}, \emptyset \} $ of $\{1, \dots, 7\}$ and consider the fixed point $p_S$ associated to the subquiver in red in Figure \ref{fig:Sing}. 
It is possible to check directly that  $\mathrm{Sing}_i (p_S)$ is non empty only if $i=3, 6$. The singular sets for $i=3,6$ are then displayed in Figure \ref{fig:SingSets}. As an immediate consequence of Formula \ref{dim}, the dimension of tangent space at $p_S$ is equal to 24. \end{example}
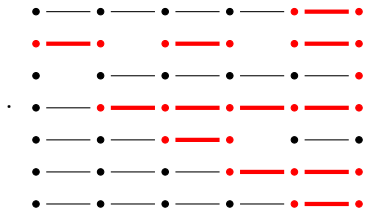
\begin{figure}[ht].
	\centering
		\begin{tikzpicture}[baseline=(current bounding box.center), scale =.85]
			\tikzstyle{point}=[circle,thick,draw=black,fill=black,inner sep=0pt,minimum width=2pt,minimum height=2pt]
			\node[] (v00) at (-1,1) {};
            \node[] (v01) at (-1,0.5) {};
			\node[] (v02) at (-1,0) {};
			\node[] (v03) at (-1,-0.5) {};
			\node[] (v04) at (-1,-1) {};
            \node[] (v05) at (-1,-1.5) {};
			\node[] (v06) at (-1,-2) {};
            \draw[fill] (v00)  circle (.05);
			\draw[fill,red] (v01)  circle (.05);
			\draw[fill] (v02)  circle (.05);
			\draw[fill] (v03)  circle (.05);
			\draw[fill] (v04)  circle (.05);
			\draw[fill] (v05)  circle (.05);
			\draw[fill] (v06)  circle (.05);
			\node[] (v10) at (0,1) {};
            \node[] (v11) at (0,0.5) {};
			\node[] (v12) at (0,0) {};
			\node[] (v13) at (0,-0.5) {};
			\node[] (v14) at (0,-1) {};
            \node[] (v15) at (0,-1.5) {};
			\node[] (v16) at (0,-2) {};
            \draw[fill] (v10)  circle (.05);
			\draw[fill,red] (v11)  circle (.05);
			\draw[fill, ] (v12)  circle (.05);
			\draw[fill,red] (v13)  circle (.05);
			\draw[fill] (v14)  circle (.05);
			\draw[fill] (v15)  circle (.05);
			\draw[fill] (v16)  circle (.05);
			\node[] (v20) at (1,1) {};
            \node[] (v21) at (1,0.5) {};
			\node[] (v22) at (1,0) {};
			\node[] (v23) at (1,-0.5) {};
			\node[] (v24) at (1,-1) {};
            \node[] (v25) at (1,-1.5) {};
			\node[] (v26) at (1,-2) {};
            \draw[fill] (v20)  circle (.05);
			\draw[fill, red] (v21)  circle (.05);
			\draw[fill] (v22)  circle (.05);
			\draw[fill,red] (v23)  circle (.05);
			\draw[fill,red] (v24)  circle (.05);
			\draw[fill] (v25)  circle (.05);
			\draw[fill] (v26)  circle (.05);
			\node[] (v30) at (2,1) {};
            \node[] (v31) at (2,0.5) {};
			\node[] (v32) at (2,0) {};
			\node[] (v33) at (2,-0.5) {};
			\node[] (v34) at (2,-1) {};
            \node[] (v35) at (2,-1.5) {};
			\node[] (v36) at (2,-2) {};
            \draw[fill] (v30)  circle (.05);
			\draw[fill,red] (v31)  circle (.05);
			\draw[fill] (v32)  circle (.05);
			\draw[fill,red] (v33)  circle (.05);
			\draw[fill,red] (v34)  circle (.05);
			\draw[fill,red] (v35)  circle (.05);
			\draw[fill] (v36)  circle (.05);
			\node[] (v40) at (3,1) {};
            \node[] (v41) at (3,0.5) {};
			\node[] (v42) at (3,0) {};
			\node[] (v43) at (3,-0.5) {};
			\node[] (v44) at (3,-1) {};
            \node[] (v45) at (3,-1.5) {};
			\node[] (v46) at (3,-2) {};
            \draw[fill,red] (v40)  circle (.05);
			\draw[fill,red] (v41)  circle (.05);
			\draw[fill] (v42)  circle (.05);
			\draw[fill,red] (v43)  circle (.05);
			\draw[fill,] (v44)  circle (.05);
			\draw[fill,red] (v45)  circle (.05);
			\draw[fill,red] (v46)  circle (.05);
			\node[] (v50) at (4,1) {};
            \node[] (v51) at (4,0.5) {};
			\node[] (v52) at (4,0) {};
			\node[] (v53) at (4,-0.5) {};
			\node[] (v54) at (4,-1) {};
            \node[] (v55) at (4,-1.5) {};
			\node[] (v56) at (4,-2) {};
            \draw[fill,red] (v50)  circle (.05);
			\draw[fill,red] (v51)  circle (.05);
			\draw[fill,red] (v52)  circle (.05);
			\draw[fill,red] (v53)  circle (.05);
			\draw[fill] (v54)  circle (.05);
			\draw[fill,red] (v55)  circle (.05);
			\draw[fill,red] (v56)  circle (.05);
			
			\draw[, ] (v00) -- (v10);
			\draw[ultra thick,red  ] (v01) -- (v11);
			\draw[,  ] (v03) -- (v13);
			\draw[, ] (v04) -- (v14);
			\draw[,  ] (v05) -- (v15);
			\draw[,  ] (v06) -- (v16);

			\draw[, ] (v10) -- (v20);
			\draw[,  ] (v12) -- (v22);
			\draw[ultra thick,red ] (v13) -- (v23);
			\draw[, ] (v14) -- (v24);
			\draw[,  ] (v15) -- (v25);
			\draw[,  ] (v16) -- (v26);
			
			\draw[, ] (v20) -- (v30);
			\draw[ultra thick,red  ] (v21) -- (v31);
			\draw[,  ] (v22) -- (v32);
			\draw[ultra thick,red  ] (v23) -- (v33);
			\draw[ultra thick,red] (v24) -- (v34);
			\draw[,  ] (v25) -- (v35);
			\draw[,  ] (v26) -- (v36);
			
			\draw[, ] (v30) -- (v40);

			\draw[,  ] (v32) -- (v42);
			\draw[ultra thick,red  ] (v33) -- (v43);

			\draw[ultra thick,red  ] (v35) -- (v45);
			\draw[,  ] (v36) -- (v46);
			
			\draw[ultra thick,red  ] (v40) -- (v50);
			\draw[ultra thick,red ] (v41) -- (v51);
			\draw[,  ] (v42) -- (v52);
			\draw[ultra thick,red  ] (v43) -- (v53);
			\draw[, ] (v44) -- (v54);
			\draw[ultra thick,red  ] (v45) -- (v55);
			\draw[ultra thick,red ] (v46) -- (v56);
		\end{tikzpicture}
			\caption{The fixed point $p_S$.}\label{fig:Sing}
\end{figure}
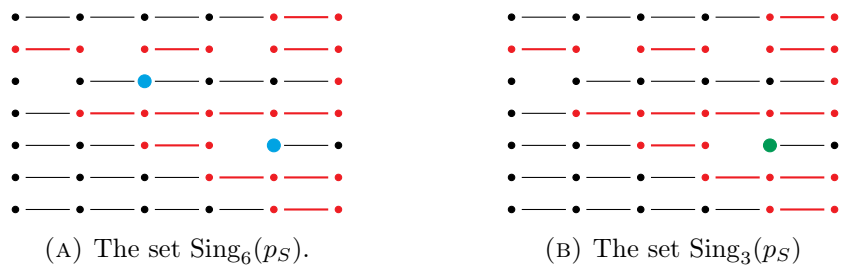
\begin{figure}[ht]
\centering
\begin{subfigure}{.3\textwidth}
\centering
		\begin{tikzpicture}[baseline=(current bounding box.center), scale =.85]
			\tikzstyle{point}=[circle,thick,draw=black,fill=black,inner sep=0pt,minimum width=2pt,minimum height=2pt]
			\node[] (v00) at (-1,1) {};
            \node[] (v01) at (-1,0.5) {};
			\node[] (v02) at (-1,0) {};
			\node[] (v03) at (-1,-0.5) {};
			\node[] (v04) at (-1,-1) {};
            \node[] (v05) at (-1,-1.5) {};
			\node[] (v06) at (-1,-2) {};
            \draw[fill] (v00)  circle (.05);
			\draw[fill,Red] (v01)  circle (.05);
			\draw[fill] (v02)  circle (.05);
			\draw[fill] (v03)  circle (.05);
			\draw[fill] (v04)  circle (.05);
			\draw[fill] (v05)  circle (.05);
			\draw[fill] (v06)  circle (.05);
			\node[] (v10) at (0,1) {};
            \node[] (v11) at (0,0.5) {};
			\node[] (v12) at (0,0) {};
			\node[] (v13) at (0,-0.5) {};
			\node[] (v14) at (0,-1) {};
            \node[] (v15) at (0,-1.5) {};
			\node[] (v16) at (0,-2) {};
            \draw[fill] (v10)  circle (.05);
			\draw[fill,Red] (v11)  circle (.05);
			\draw[fill,] (v12)  circle (.05);
			\draw[fill,Red] (v13)  circle (.05);
			\draw[fill] (v14)  circle (.05);
			\draw[fill] (v15)  circle (.05);
			\draw[fill] (v16)  circle (.05);
			\node[] (v20) at (1,1) {};
            \node[] (v21) at (1,0.5) {};
			\node[] (v22) at (1,0) {};
			\node[] (v23) at (1,-0.5) {};
			\node[] (v24) at (1,-1) {};
            \node[] (v25) at (1,-1.5) {};
			\node[] (v26) at (1,-2) {};
            \draw[fill] (v20)  circle (.05);
			\draw[fill,Red] (v21)  circle (.05);
			\draw[fill, Cerulean] (v22)  circle (.1);
			\draw[fill,Red] (v23)  circle (.05);
			\draw[fill,Red] (v24)  circle (.05);
			\draw[fill] (v25)  circle (.05);
			\draw[fill] (v26)  circle (.05);
			\node[] (v30) at (2,1) {};
            \node[] (v31) at (2,0.5) {};
			\node[] (v32) at (2,0) {};
			\node[] (v33) at (2,-0.5) {};
			\node[] (v34) at (2,-1) {};
            \node[] (v35) at (2,-1.5) {};
			\node[] (v36) at (2,-2) {};
            \draw[fill] (v30)  circle (.05);
			\draw[fill,Red] (v31)  circle (.05);
			\draw[fill] (v32)  circle (.05);
			\draw[fill,Red] (v33)  circle (.05);
			\draw[fill,Red] (v34)  circle (.05);
			\draw[fill,Red] (v35)  circle (.05);
			\draw[fill] (v36)  circle (.05);
			\node[] (v40) at (3,1) {};
            \node[] (v41) at (3,0.5) {};
			\node[] (v42) at (3,0) {};
			\node[] (v43) at (3,-0.5) {};
			\node[] (v44) at (3,-1) {};
            \node[] (v45) at (3,-1.5) {};
			\node[] (v46) at (3,-2) {};
            \draw[fill,Red] (v40)  circle (.05);
			\draw[fill,Red] (v41)  circle (.05);
			\draw[fill] (v42)  circle (.05);
			\draw[fill,Red] (v43)  circle (.05);
			\draw[fill, Cerulean] (v44)  circle (.1);
			\draw[fill,Red] (v45)  circle (.05);
			\draw[fill,Red] (v46)  circle (.05);
			\node[] (v50) at (4,1) {};
            \node[] (v51) at (4,0.5) {};
			\node[] (v52) at (4,0) {};
			\node[] (v53) at (4,-0.5) {};
			\node[] (v54) at (4,-1) {};
            \node[] (v55) at (4,-1.5) {};
			\node[] (v56) at (4,-2) {};
            \draw[fill,Red] (v50)  circle (.05);
			\draw[fill,Red] (v51)  circle (.05);
			\draw[fill,Red] (v52)  circle (.05);
			\draw[fill,Red] (v53)  circle (.05);
			\draw[fill] (v54)  circle (.05);
			\draw[fill,Red] (v55)  circle (.05);
			\draw[fill,Red] (v56)  circle (.05);
			
			\draw[, ] (v00) -- (v10);
			\draw[thick,Red  ] (v01) -- (v11);
			\draw[,  ] (v03) -- (v13);
			\draw[, ] (v04) -- (v14);
			\draw[,  ] (v05) -- (v15);
			\draw[,  ] (v06) -- (v16);

			\draw[, ] (v10) -- (v20);
			\draw[,  ] (v12) -- (v22);
			\draw[thick,Red ] (v13) -- (v23);
			\draw[, ] (v14) -- (v24);
			\draw[,  ] (v15) -- (v25);
			\draw[,  ] (v16) -- (v26);
			
			\draw[, ] (v20) -- (v30);
			\draw[thick,Red  ] (v21) -- (v31);
			\draw[,  ] (v22) -- (v32);
			\draw[thick,Red  ] (v23) -- (v33);
			\draw[thick, Red] (v24) -- (v34);
			\draw[,  ] (v25) -- (v35);
			\draw[,  ] (v26) -- (v36);
			
			\draw[, ] (v30) -- (v40);

			\draw[,  ] (v32) -- (v42);
			\draw[thick,Red  ] (v33) -- (v43);

			\draw[thick,Red  ] (v35) -- (v45);
			\draw[,  ] (v36) -- (v46);
			
			\draw[thick,Red  ] (v40) -- (v50);
			\draw[thick,Red ] (v41) -- (v51);
			\draw[,  ] (v42) -- (v52);
			\draw[thick,Red  ] (v43) -- (v53);
			\draw[, ] (v44) -- (v54);
			\draw[thick,Red  ] (v45) -- (v55);
			\draw[thick,Red ] (v46) -- (v56);
		\end{tikzpicture}
			\caption{The set $\mathrm{Sing}_6(p_S)$.}
\end{subfigure}\hspace{0.1\textwidth}
\begin{subfigure}{.3\textwidth}
\centering
\begin{tikzpicture}[baseline=(current bounding box.center), scale =.85]
			\tikzstyle{point}=[circle,thick,draw=black,fill=black,inner sep=0pt,minimum width=2pt,minimum height=2pt]
			\node[] (v00) at (-1,1) {};
            \node[] (v01) at (-1,0.5) {};
			\node[] (v02) at (-1,0) {};
			\node[] (v03) at (-1,-0.5) {};
			\node[] (v04) at (-1,-1) {};
            \node[] (v05) at (-1,-1.5) {};
			\node[] (v06) at (-1,-2) {};
            \draw[fill] (v00)  circle (.05);
			\draw[fill,Red] (v01)  circle (.05);
			\draw[fill] (v02)  circle (.05);
			\draw[fill] (v03)  circle (.05);
			\draw[fill] (v04)  circle (.05);
			\draw[fill] (v05)  circle (.05);
			\draw[fill] (v06)  circle (.05);
			\node[] (v10) at (0,1) {};
            \node[] (v11) at (0,0.5) {};
			\node[] (v12) at (0,0) {};
			\node[] (v13) at (0,-0.5) {};
			\node[] (v14) at (0,-1) {};
            \node[] (v15) at (0,-1.5) {};
			\node[] (v16) at (0,-2) {};
            \draw[fill] (v10)  circle (.05);
			\draw[fill,Red] (v11)  circle (.05);
			\draw[fill, ] (v12)  circle (.05);
			\draw[fill,Red] (v13)  circle (.05);
			\draw[fill] (v14)  circle (.05);
			\draw[fill] (v15)  circle (.05);
			\draw[fill] (v16)  circle (.05);
			\node[] (v20) at (1,1) {};
            \node[] (v21) at (1,0.5) {};
			\node[] (v22) at (1,0) {};
			\node[] (v23) at (1,-0.5) {};
			\node[] (v24) at (1,-1) {};
            \node[] (v25) at (1,-1.5) {};
			\node[] (v26) at (1,-2) {};
            \draw[fill] (v20)  circle (.05);
			\draw[fill,Red] (v21)  circle (.05);
			\draw[fill] (v22)  circle (.05);
			\draw[fill,Red] (v23)  circle (.05);
			\draw[fill,Red] (v24)  circle (.05);
			\draw[fill] (v25)  circle (.05);
			\draw[fill] (v26)  circle (.05);
			\node[] (v30) at (2,1) {};
            \node[] (v31) at (2,0.5) {};
			\node[] (v32) at (2,0) {};
			\node[] (v33) at (2,-0.5) {};
			\node[] (v34) at (2,-1) {};
            \node[] (v35) at (2,-1.5) {};
			\node[] (v36) at (2,-2) {};
            \draw[fill] (v30)  circle (.05);
			\draw[fill,Red] (v31)  circle (.05);
			\draw[fill] (v32)  circle (.05);
			\draw[fill,Red] (v33)  circle (.05);
			\draw[fill,Red] (v34)  circle (.05);
			\draw[fill,Red] (v35)  circle (.05);
			\draw[fill] (v36)  circle (.05);
			\node[] (v40) at (3,1) {};
            \node[] (v41) at (3,0.5) {};
			\node[] (v42) at (3,0) {};
			\node[] (v43) at (3,-0.5) {};
			\node[] (v44) at (3,-1) {};
            \node[] (v45) at (3,-1.5) {};
			\node[] (v46) at (3,-2) {};
            \draw[fill,Red] (v40)  circle (.05);
			\draw[fill,Red] (v41)  circle (.05);
			\draw[fill] (v42)  circle (.05);
			\draw[fill,Red] (v43)  circle (.05);
			\draw[fill, ForestGreen] (v44)  circle (.1);
			\draw[fill,Red] (v45)  circle (.05);
			\draw[fill,Red] (v46)  circle (.05);
			\node[] (v50) at (4,1) {};
            \node[] (v51) at (4,0.5) {};
			\node[] (v52) at (4,0) {};
			\node[] (v53) at (4,-0.5) {};
			\node[] (v54) at (4,-1) {};
            \node[] (v55) at (4,-1.5) {};
			\node[] (v56) at (4,-2) {};
            \draw[fill,Red] (v50)  circle (.05);
			\draw[fill,Red] (v51)  circle (.05);
			\draw[fill,Red] (v52)  circle (.05);
			\draw[fill,Red] (v53)  circle (.05);
			\draw[fill] (v54)  circle (.05);
			\draw[fill,Red] (v55)  circle (.05);
			\draw[fill,Red] (v56)  circle (.05);
			
			\draw[, ] (v00) -- (v10);
			\draw[thick,Red  ] (v01) -- (v11);
			\draw[,  ] (v03) -- (v13);
			\draw[, ] (v04) -- (v14);
			\draw[,  ] (v05) -- (v15);
			\draw[,  ] (v06) -- (v16);

			\draw[, ] (v10) -- (v20);
			\draw[,  ] (v12) -- (v22);
			\draw[thick,Red ] (v13) -- (v23);
			\draw[, ] (v14) -- (v24);
			\draw[,  ] (v15) -- (v25);
			\draw[,  ] (v16) -- (v26);
			
			\draw[, ] (v20) -- (v30);
			\draw[thick,Red  ] (v21) -- (v31);
			\draw[,  ] (v22) -- (v32);
			\draw[thick,Red  ] (v23) -- (v33);
			\draw[thick,Red] (v24) -- (v34);
			\draw[,  ] (v25) -- (v35);
			\draw[,  ] (v26) -- (v36);
			
			\draw[, ] (v30) -- (v40);

			\draw[,  ] (v32) -- (v42);
			\draw[thick,Red  ] (v33) -- (v43);

			\draw[thick,Red  ] (v35) -- (v45);
			\draw[,  ] (v36) -- (v46);
			
			\draw[thick,Red  ] (v40) -- (v50);
			\draw[thick,Red ] (v41) -- (v51);
			\draw[,  ] (v42) -- (v52);
			\draw[thick,Red  ] (v43) -- (v53);
			\draw[, ] (v44) -- (v54);
			\draw[thick,Red  ] (v45) -- (v55);
			\draw[thick,Red ] (v46) -- (v56);
		\end{tikzpicture}
			\caption{The set $\mathrm{Sing}_3(p_S)$}
\end{subfigure}
\caption{The singularity sets of point $p_S$.}\label{fig:SingSets}
\end{figure}
\section{Criteria for Smoothness at Fixed Points}\label{sec:4}
In this section we provide a list of equivalent smoothness criteria for a fixed point $p_S$. More closely, we identify a property of the sequence $S$ that generalizes the smoothness condition proved in \cite{CII} for Feigin degenerations. Moreover to each $p_S$ is associated a graph encoding informations about dimension of tangent space at $p_S$. Our results are proved in Section~\ref{sect:proofs}
\subsection{The Generalized CI-F-R Condition.}\label{subsec:GenCIFR}
\begin{figure}[ht]
\centering
\begin{subfigure}{.3\textwidth}
\centering
		\begin{tikzpicture}[baseline=(current bounding box.center), scale =.8]
			\tikzstyle{point}=[circle,thick,draw=black,fill=black,inner sep=0pt,minimum width=2pt,minimum height=2pt]
			\node[] (v00) at (-1,1) {};
            \node[] (v01) at (-1,0) {};
			\node[] (v02) at (-1,-1) {};
			\draw[fill] (v00)  circle (.07);
			\draw[fill] (v01)  circle (.07);
			\draw[fill] (v02)  circle (.07);
			\node[] (v10) at (1,1) {};
			\node[] (v11) at (1,0) {};
			\node[] (v12) at (1,-1) {};
			\draw[fill] (v10)  circle (.07);
			\draw[fill] (v11)  circle (.07);
			\draw[fill] (v12)  circle (.07);
		\end{tikzpicture}
\end{subfigure}\hspace{0.05\textwidth}
\begin{subfigure}{.3\textwidth}
\centering
		\begin{tikzpicture}[baseline=(current bounding box.center), scale =.8]
			\tikzstyle{point}=[circle,thick,draw=black,fill=black,inner sep=0pt,minimum width=2pt,minimum height=2pt]
			\node[] (v00) at (-1,1) {};
            \node[] (v01) at (-1,0) {};
			\node[] (v02) at (-1,-1) {};
			\draw[fill] (v00)  circle (.07);
			\draw[fill] (v01)  circle (.07);
			\draw[fill,red] (v02)  circle (.07);
			\node[] (v10) at (1,1) {};
			\node[] (v11) at (1,0) {};
			\node[] (v12) at (1,-1) {};
			\draw[fill] (v10)  circle (.07);
			\draw[fill,red] (v11)  circle (.07);
			\draw[fill,red] (v12)  circle (.07);
		\end{tikzpicture}
\end{subfigure}
\caption{The coefficient quiver of $\mathbb{P}^2(\C) \times \mathbb{P}^2 (\C)$ and the fixed point associated to  $S=\{\{1\}, \{1,2\}\}$.}\label{fig:P2P2}
\end{figure}
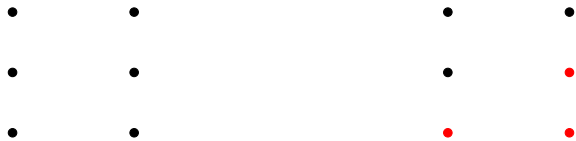

We recall that Feigin degeneration is the ${\bf R}$ degeneration of $ {\mathcal{F}l}(\C^{n+1})$ such that for every $i$ we have $R_i=\{i+1\}$. A fixed point $p_S$ for the torus action described in Section~\ref{sec:3} is associated to an admissible sequence $S=(S_1, \dots, S_n)$ such that $S_i \subset S_{i+1} \cup \{i+1\}$. In \cite{CII} the fixed points in the smooth locus are classified and enumerated, providing the following combinatorial smoothness criterion.
\begin{defn}[Cerulli Irelli - Feigin - Reineke Condition]\label{CIFR}
We say that an admissible sequence $S$ for the Feigin degeneration has the Cerulli Irelli - Feigin - Reineke Condition (for short CI-F-R Condition) if for all $h,k$ such that $1 \leq h<k \leq n$ we have:
\begin{equation} k \in S_{h} \Rightarrow h+1 \in S_{k}. \end{equation}
\end{defn}
\begin{crit}[Smoothness Criterion - Feigin Degeneration, \cite{CII}, Theorem 4.2]\label{crit:CIFR}
A fixed point $p_S$ is smooth if and only if $S$ has the CI-F-R condition.
\end{crit}
Our first smoothness criterion is the natural generalization of this one.
\begin{defn}[Generalized CI-F-R Condition]\label{GenCIFR}
Let $S$ be an $\bf{R}$-admissible sequence of subsets of $\{1, \dots, n+1\}$. We say that $S$ has the Generalized CI-F-R Condition if, for every $i \in \{1, \dots, n+1\}$, 
\begin{equation} i \in S_{k}, \, k \in \pos (j) \Rightarrow j \in S_{h+1}, \; \forall h \in \pos(i), \, h\geq k.\end{equation}
\end{defn}
Observe that in the Feigin degeneration we have $\pos (i)=\{i-1\}$ for all $i$. Consequently an admissible sequence $S$ has Generalized CI-F-R Condition if for all $j,h$ such that $j \leq h$ we have  $h \in S_j \Rightarrow j+1 \in S_h$, that agrees with the CI-F-R condition.
\begin{proposition}\label{CIFREquivalence}
An $\bf{R}$-admissible sequence $S$ has the Generalized CI-F-R Condition if and only if the set $ \mathrm{Sing} (p_S)$ is empty. 
\end{proposition}
\proof It is an immediate consequence of the definition that if $S$ has the Generalized CI-F-R Condition then $\mathrm{Sing} (p_S)=\emptyset$.
On the other hand, $\mathrm{Sing} (p_S)=\emptyset$ is equivalent to requiring that  $\mathrm{Sing}_i (p_S)=\emptyset$ for all $i$. Suppose now that there exist $i, j$ and $ k \in \pos(j)$, $h \in \pos(i)$ such that $k \leq h$ and $i \in S_{k}, j \notin S_{h+1}$. For a fixed pair $i, j$, consider $h$ minimal with this property. First of all, observe that $j$ must be different from $i$, otherwise $(i,h+1)$ is in $\mathrm{Sing}_i (p_S)$. Note that in particular this implies that there exists a minimal index $k_i$ such that $i \in S_t$ for all $t \geq k_i$. So, we are in the case $i \neq j$ and set $\widetilde{k}=\max\{t \in \pos(j)\, | \, t \leq h\}$. The segment $[\widetilde{k}+1,h+1]^j_M$ is connected by definition of $\widetilde{k}$ and consequently $j \notin S_t$ for all $t\in \{\widetilde{k}+1, \dots, h+1\}$. By minimality of $h$ this implies that the segment $[\widetilde{k},h]^i_S$ is connected too. Moreover, because $k \leq  \widetilde{k}$,  the fact that $i \in S_t$ for all $t \geq k_i$ implies that $i \in S_t, \; \forall t \in \{\widetilde{k}, \dots, h\}$ and $(j,h+1)$ satisfies the conditions of Definition \ref{SingSet}, i.e. $\mathrm{Sing}_i (p_S) \neq \emptyset$, that is absurd. \endproof
As a consequence of Formula~\ref{dim}, if ${\mathcal{F}l}^\mathbf{R}(\C^{n+1})$ is a flat degeneration, one obtains the following smoothness criterion that generalize Criterion~\ref{crit:CIFR}:
\begin{corollary}\label{GenCifr} Let ${\mathcal{F}l}^\mathbf{R}(\C^{n+1})$ be a flat degeneration. A fixed point $p_S$ is smooth if and only if $S$ satisfies the Generalized CI-F-R Condition.
\end{corollary}
\begin{rmk}\label{saturation}
Suppose that $S$ has the Generalized CI-F-R Condition, then $i \in S_h$ implies $i \in S_k$ for all $k$ greater or equal than $h$, or equivalently that $S_i \subset S_{i+1}$ for every $i$. The converse is not true, in fact the admissible sequence $S=\{\{3\},\{3,4\},\{1,3,4\}\}$ defining the successor closed subquiver in Figure~\ref{fig:sing3} satisfies the condition $S_i \subset S_{i+1}$ for every $i$, but as a consequence of Formula~(\ref{dim}) and of Example~\ref{ex:nonsmoothF3}, the associated fixed point $p_S$ is not  smooth.
\end{rmk}
\begin{example}\label{ex:P2P2}
The next example shows that Corollary \ref{GenCifr} does not hold if the degeneration is not flat. The algebraic variety
$\mathbb{P}^2(\C) \times \mathbb{P}^2 (\C)$ can be identified by the degeneration ${\mathcal{F}l}^{\bf R}(\C^3)$ such that ${\bf R}=(\{1,2,3\})$. In particular, it is not a flat degeneration. In Figure \ref{fig:P2P2} the coefficient quiver of $\mathbb{P}^2(\C) \times \mathbb{P}^2 (\C)$ is displayed. Moreover, the second diagram is the successor closed subquiver associated to the point $p_S$ defined by the admissible sequence $S=(\{1\}, \{1,2\})$. The set of singularities of $p_S$ is non empty, in fact $\mathrm{Sing}(p_S)=\mathrm{Sing}_1(p_S)=\{(3,2)\}$, and then $S$ does not have the Generalized CI-F-R Condition. However, 
$\mathbb{P}^2(\C) \times \mathbb{P}^2 (\C)$ is a smooth algebraic variety and so every fixed point is smooth. Nevertheless, we remark that if we compute the dimension of tangent space at $p_S$ using Equation~(\ref{dim}), our formula produces the correct result. \end{example}
\subsection{The Mutation Graph}\label{subsec:graph}
We introduce now a graph that encodes some informations about the mutations that occur between the coefficient quiver of a fixed point and its adjacent vertices in the moment graph. 
We prove that the topology of this graph is linked with the property of a fixed point $p_S$ of being smooth. \\
We say that we have a  mutation between the $i$-th row and the $j$-th row of $Q_S$  if there exists a fixed point $S'$ such that $Q_{S'}$ is obtained from $Q_S$ moving a unique movable part from its $i$-th row to the $j$-th row. 
For a fixed subquiver $Q_S$, it is possible to have many mutations between row $i$ and row $j$. We will denote the set of these mutations by $\mathrm{Mut}_S(i,j)$. 
\begin{defn}[Oriented Mutation Graphs]\label{defn:mutG}
Consider a fixed point $p_S \in {\mathcal{F}l}^\mathbf{R}(\C^{n+1})$. The oriented mutation graph of $p_S$ is the oriented multigraph $\widetilde{G}_S=(V(\widetilde{G}_S),E(\widetilde{G}_S))$ described by the following data:
\begin{enumerate}
    \item The set of vertices $V(\widetilde{G}_S)$ is the set $\{1, \dots , n+1\}$,
    \item There is an oriented edge from $i$ to $j$ for every mutation in $\mathrm{Mut}_S(i,j)$.
\end{enumerate}
\end{defn}
We are going to extensively use the undirected version $G_S$ of the graph $\widetilde{G}_S$, i.e. the same graph where we forgot the edges orientations. We will refer to $G_S$ as the \emph{Mutation Graph} of $S$. 
\begin{example}In Figure \ref{fig:fixpointmut} a fixed point for the Feigin Degeneration when $n=3$ is displayed.
\begin{figure}[ht]
	\centering
		\begin{tikzpicture}[baseline=(current bounding box.center), scale =.85]
			\tikzstyle{point}=[circle,thick,draw=black,fill=black,inner sep=0pt,minimum width=2pt,minimum height=2pt]
			\node[] (v00) at (-1,1) {};
            \node[] (v01) at (-1,0.333) {};
			\node[] (v02) at (-1,-0.330) {};
			\node[] (v03) at (-1,-1) {};
			\draw[fill,red] (v00)  circle (.07);
			\draw[fill] (v01)  circle (.05);
			\draw[fill] (v02)  circle (.05);
			\draw[fill] (v03)  circle (.05);
			\node[] (v10) at (1,1) {};
			\node[] (v11) at (1,0.333) {};
			\node[] (v12) at (1,-0.330) {};
			\node[] (v13) at (1,-1) {};
			\draw[fill,red] (v10)  circle (.07);
			\draw[fill,red] (v11)  circle (.07);
			\draw[fill] (v12)  circle (.05);
			\draw[fill] (v13)  circle (.05);
			\node[] (v20) at (3,1) {};
			\node[] (v21) at (3,0.333) {};
			\node[] (v22) at (3,-0.330) {};
			\node[] (v23) at (3,-1) {};
			\draw[fill,red] (v20)  circle (.07);
			\draw[fill,red] (v21)  circle (.07);
			\draw[fill,red] (v22)  circle (.07);
			\draw[fill] (v23)  circle (.05);
			
			\draw[ultra thick, red] (v00) -- (v10);
			\draw[,  ] (v01) -- (v11);
			\draw[,  ] (v03) -- (v13);

			\draw[ultra thick, red] (v10) -- (v20);
			\draw[,  ] (v12) -- (v22);
			\draw[,  ] (v13) -- (v23);
		\end{tikzpicture}
			\caption{The fixed point $p_S$ defined by $S=(\{4\}, \{3,4\}, \{2,3,4\})$}\label{fig:fixpointmut}
\end{figure}
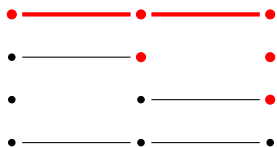
The associated oriented mutation graph and its unoriented version are displayed in Figure \ref{fig:mutgraph}.
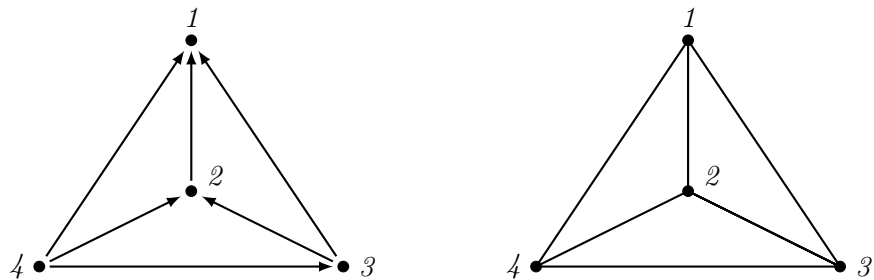
\begin{figure}[ht]
\centering
\begin{subfigure}{.3\textwidth}
\centering
\begin{tikzpicture}
\node (a) at (0,0) {};
\node (b) at (4,0) {};
\node (c) at (2,1){};
\node (d) at (2,3){};

\draw[fill=black] (0,0) circle (2pt);
\draw[fill=black] (4,0) circle (2pt);
\draw[fill=black] (2,1) circle (2pt);
\draw[fill=black] (2,3) circle (2pt);
\node at (-0.3,0) {4};
\node at (4.3,0) {3};
\node at (2.3,1.2) {2};
\node at (2,3.3) {1};
    \draw[thick,-latex] (a) -- (b);
    \draw[thick,-latex](a) -- (c);
    \draw[thick,-latex] (a) -- (d);
    \draw[thick,-latex] (b) -- (c);
    \draw[thick,-latex](b) -- (d);
    \draw[thick,-latex] (c) -- (d);
\end{tikzpicture}
\end{subfigure}\hspace{0.1\textwidth}
\begin{subfigure}{.3\textwidth}
\centering
\begin{tikzpicture}
\draw[fill=black] (0,0) circle (2pt);
\draw[fill=black] (4,0) circle (2pt);
\draw[fill=black] (2,1) circle (2pt);
\draw[fill=black] (2,3) circle (2pt);
\node at (-0.3,0) {4};
\node at (4.3,0) {3};
\node at (2.3,1.2) {2};
\node at (2,3.3) {1};
\draw[thick](0,0) -- (4,0) -- (2,1) -- (0,0) -- (2,3) -- (4,0) -- (2,1) -- (2,3);
\end{tikzpicture}
\end{subfigure}
\caption{The oriented and unoriented mutation graphs of $p_S$ in Figure \ref{fig:fixpointmut}}\label{fig:mutgraph}
\end{figure}
\end{example}
\begin{example}
In this example we exhibit a fixed point for a flat degeneration that has a non simple mutation graph. In Figure~\ref{fig:nonsimplflat} is displayed the successor closed subquiver associated to the admissible sequence $S=\{\{3\},\{3,4\},\{1,3,4\}\}$, that defines a $T$-fixed point $p_S$ for the Feigin degeneration of ${\mathcal{F}l}(\C^4)$. It is immediate to check that the set $\mrm{Mut}(3,2)$ has 2 elements. The associated oriented mutation graph is a multigraph and it is displayed on the right side of Figure~\ref{fig:nonsimplflat}. \end{example}
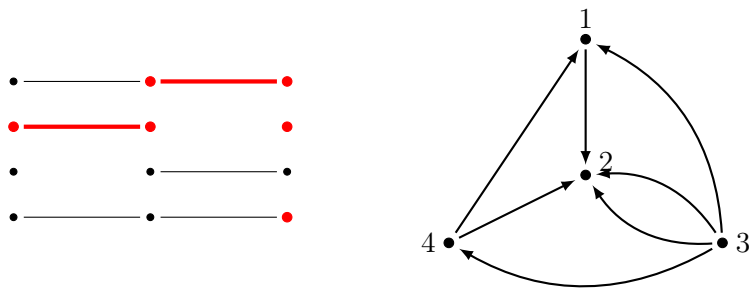
\begin{figure}[ht]
\centering
\begin{subfigure}{.3\textwidth}
\centering
			\begin{tikzpicture}[baseline=(current bounding box.center), scale =.9]
			\tikzstyle{point}=[circle,thick,draw=black,fill=black,inner sep=0pt,minimum width=2pt,minimum height=2pt]
			\node[] (v00) at (-1,1) {};
            \node[] (v01) at (-1,0.333) {};
			\node[] (v02) at (-1,-0.330) {};
			\node[] (v03) at (-1,-1) {};
			\draw[fill] (v00)  circle (.05);
			\draw[fill,red] (v01)  circle (.07);
			\draw[fill] (v02)  circle (.05);
			\draw[fill] (v03)  circle (.05);
			\node[] (v10) at (1,1) {};
			\node[] (v11) at (1,0.333) {};
			\node[] (v12) at (1,-0.330) {};
			\node[] (v13) at (1,-1) {};
			\draw[fill,red] (v10)  circle (.07);
			\draw[fill,red] (v11)  circle (.07);
			\draw[fill] (v12)  circle (.05);
			\draw[fill] (v13)  circle (.05);
			\node[] (v20) at (3,1) {};
			\node[] (v21) at (3,0.333) {};
			\node[] (v22) at (3,-0.330) {};
			\node[] (v23) at (3,-1) {};
			\draw[fill,red] (v20)  circle (.07);
			\draw[fill,red] (v21)  circle (.07);
			\draw[fill] (v22)  circle (.05);
			\draw[fill,red] (v23)  circle (.07);

			\draw[] (v00) -- (v10);
			\draw[ultra thick, red] (v01) -- (v11);
			\draw[,  ] (v03) -- (v13);

			\draw[ultra thick, red] (v10) -- (v20);
			\draw[,  ] (v12) -- (v22);
			\draw[,  ] (v13) -- (v23);
		\end{tikzpicture}
\end{subfigure}\hspace{0.05\textwidth}
\begin{subfigure}{.3\textwidth}
\centering
\begin{tikzpicture}[baseline=(current bounding box.center), scale =.9]
\node (a) at (0,0) {};
\node (b) at (4,0) {};
\node (c) at (2,1){};
\node (d) at (2,3){};

\draw[fill=black] (0,0) circle (2pt);
\draw[fill=black] (4,0) circle (2pt);
\draw[fill=black] (2,1) circle (2pt);
\draw[fill=black] (2,3) circle (2pt);
\node at (-0.3,0) {4};
\node at (4.3,0) {3};
\node at (2.3,1.2) {2};
\node at (2,3.3) {1};
    \draw[thick,-latex] (b) to [bend left] (a);
    \draw[thick,-latex](a) -- (c);
    \draw[thick,-latex] (a) -- (d);
    \draw[thick,-latex](b) to [bend right] (d);
    \draw[thick,-latex] (d) -- (c);
\draw [thick, -latex] (b) to [bend right] (c);
\draw [thick, -latex] (b) to [bend left] (c);
\end{tikzpicture}
\end{subfigure}
\caption{A fixed point $p_S$ for Feigin degeneration with non simple oriented mutation graph.}\label{fig:nonsimplflat}
\end{figure}
A consequence of our work is that certain combinatorial properties of $G_S$ and $\widetilde{G}_S$ are equivalent to the fact that $p_S$ is smooth. 
Nevertheless, an explicit link between the structure of $G_S$ and $\widetilde{G}_S$ and geometric properties of $p_S$ is far to be clear. 
As an example, it results that all smooth fixed points for a fixed degeneration ${\mathcal{F}l}^\mathbf{R}(\C^{n+1})$  have the same oriented mutation graphs up to vertices relabeling. Moreover, in the case of complete flag variety ${\mathcal{F}l}(\C^{n+1})$ we observed that starting from the mutation graph of a fixed point $p_S$ it is possible to recover the inversion set of the permutation $w \in S_{n+1}$ associated to $p_S$ and then the dimension of the corresponding cell. These phenomena are quite interesting and we hope to investigate deeper links between topology of mutation graphs and cell decomposition
in some future works.
\subsection{Smoothness Criteria}
\begin{figure}[ht]
\centering
\begin{subfigure}{.3\textwidth}
\centering
		\begin{tikzpicture}[baseline=(current bounding box.center), scale =.8]
			\tikzstyle{point}=[circle,thick,draw=black,fill=black,inner sep=0pt,minimum width=2pt,minimum height=2pt]
			\node[] (v00) at (-1,1) {};
            \node[] (v01) at (-1,0) {};
			\node[] (v02) at (-1,-1) {};
			\draw[fill] (v00)  circle (.07);
			\draw[fill] (v01)  circle (.07);
			\draw[fill,red] (v02)  circle (.07);
			\node[] (v10) at (1,1) {};
			\node[] (v11) at (1,0) {};
			\node[] (v12) at (1,-1) {};
			\draw[fill] (v10)  circle (.07);
			\draw[fill,red] (v11)  circle (.07);
			\draw[fill,red] (v12)  circle (.07);
		\end{tikzpicture}
\end{subfigure}\hspace{0.05\textwidth}
\begin{subfigure}{.3\textwidth}
\centering
\begin{tikzpicture}[baseline=(current bounding box.center), scale =.8]
\node (a) at (0,0) {};
\node (b) at (4,0) {};
\node (d) at (2,3){};

\draw[fill=black] (0,0) circle (2pt);
\draw[fill=black] (4,0) circle (2pt);
\draw[fill=black] (2,3) circle (2pt);
\node at (-0.3,0) {1};
\node at (4.3,0) {3};
\node at (2,3.3) {2};
    \draw[thick,-latex] (a) to [bend left] (b);
       \draw[thick,-latex] (a) to [bend right] (b);
    \draw[thick,-latex] (a) -- (d);
    \draw[thick,-latex](d) to (b);
\end{tikzpicture}
\end{subfigure}
\caption{A smooth point $p_S$ with non simple mutation graph.}\label{fig:nonsimplemut}
\end{figure}
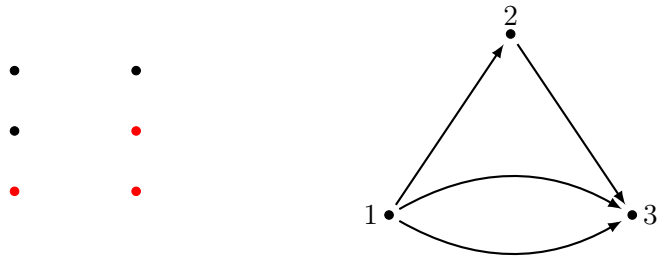

Using the tools introduced in Section \ref{subsec:GenCIFR} and \ref{subsec:graph} it is possible to provide a list of smoothness criteria for a $T$-fixed point $p_S$, summarized in Theorem~\ref{thm:criteria}.  
\begin{defn}
A \emph{Tournament Graph} $T$ over $n$ vertices is the assignment of a direction of the edges of the complete graph over $n$ vertices. We say that $T$ is \emph{transitive} if the edge directions induce a total ordering on the vertices.
\end{defn}
\begin{theorem}[Smoothness Criteria]\label{thm:criteria} Let $p_S$ be a $T$ fixed point for a flat degeneration ${\mathcal{F}l}^\mathbf{R}(\C^{n+1})$ associated to an $\bf{R}$-admissible sequence $S$. The following conditions are equivalent:
\begin{enumerate}
    \item The point $p_S$ is smooth;
    \item The set $\mrm{Sing}(p_S)$ is empty;
        \item The $\bf{R}$-admissible sequence $S$ has the Generalized CI-F-R Condition.
         \item The Mutation graph $G_S$ is the complete graph over $n+1$ vertices.
\item The Oriented Mutation graph $\widetilde{G}_S$ is a transitive tournament $n+1$ vertices;
\end{enumerate}
\end{theorem}
The equivalence of conditions $(1)$ and $(2)$ is a direct consequence of Formula~(\ref{dim}). The fact that $(2) \Longleftrightarrow (3)$ follows by Proposition~\ref{CIFREquivalence}.
Moreover the fact that $(5) \Rightarrow (4) \Rightarrow (1)$ is a consequence of flatness, of Proposition~\ref{ValSmooth} and of definition of transitive tournaments. In Section~\ref{subsec:proofgraph} we give a direct proof of the fact that $(1) \Rightarrow (5)$ using techniques coming from representation theory of quivers of type $A$. 
\begin{example}
We provide here an example of a smooth fixed point in a non flat degeneration with non a simple mutation graph. We consider the fixed point $p_S$ of Example~\ref{ex:P2P2}. We have just proved that $p_S$ is smooth. Moreover it is easy to check that $|\mrm{Mut}(1,3)|=2$ and so its mutation graph is not simple. In Figure~\ref{fig:nonsimplemut} it is displayed the successor closed subquiver associated to $p_S$ and its oriented mutation graph. We finally remark that, up to vertices relabeling, this graph is the same for every $T$-fixed point in~$\mathbb{P}^2(\C) \times \mathbb{P}^2(\C)$. \end{example}
\section{Proof of our smoothness criteria}\label{sect:proofs}
In this section we complete the proof of equivalences in the statement of Theorem~\ref{thm:criteria}. Our proofs work by induction and use extensively the combinatorics of coefficient quivers and of their successor closed subquivers. In Section~\ref{Restriction} we develop some techniques that allow us to restrict to suitable subquivers and then proceed by induction.
\subsection{Restriction of Coefficient Quiver}\label{Restriction}
Let $p_S$ be a $T$ fixed point in a linear degeneration ${\mathcal{F}l}^\mathbf{R}(\C^{n+1})$ associated to an $\bf{R}$-admissible sequence $S$. Moreover, let $M^{\bf R}$ be the quiver representation of the equioriented quiver $A_n$ associated to the degeneration ${\mathcal{F}l}^\mathbf{R}(\C^{n+1})$.
\begin{defn}
We say that $Q_S$ \emph{has an empty row} if there exist $j\in \{1, \dots, n+1\}$ such that $j \notin S_i$ for all $i$.
\end{defn}
Observe that this is coherent with the fact that, if we see $Q_S$ as a subquiver of $Q(M^{\bf R})$, no connected components of $Q_S$ belongs to the $j$-th row of $Q(M^{\bf R})$. Our induction is based on the possibility, if $Q_S$ has an empty row, of restrict mutations of $Q_S$ to mutations of a certain subquiver of $Q(M^{\bf R})$, that corresponds to a fixed point for a torus action on suitable quiver Grassmanian $\mrm{Gr}_{d}(M^{\bf R'})$, where $M^{\bf R'}$ is a representation of the equioriented quiver of type $A_{n-1}$. We now describe our inductive reasoning in a more accurate way. \\
Suppose that the $j$-th row of the quiver $Q_S$ is empty and consider the family of sets ${\bf{R'}}=(R'_1, \dots R'_{n-2})$ such that 
\begin{equation*}
R'_i=
\begin{cases}
R_i \setminus \{j\} & \text{if} \;j \in R_i \\
R_i & \text{otherwise} 
\end{cases}
\end{equation*}
The family $\bf{R'}$ defines a linear degeneration ${\mathcal{F}l}^{\bf{R'}}(\C^n)$ of the flag variety ${\mathcal{F}l}(\C^n)$. Up to identifying $\C^n$ as the subspace of $\C^{n+1}$ spanned by the vectors $\{v_h^i\}_{i \neq j}$ for each $h$, the coefficient quiver $Q(M^{\bf R'})$ associated ${\mathcal{F}l}^{\bf R'}(\C^n)$ is the subquiver of $Q(M^{\bf R})$ obtained deleting the $j$-th row and the column corresponding to the vectors $\{v_n^i\}$. Since the $j$-th row of $Q_S$ is empty, starting from the $\bf{R}$-admissible sequence $S$ it is possible to define the $\bf{R'}$-admissible sequence $S'$ setting  $S'_i=S_i $ for all~$i\leq n-1$.
\begin{rmk}\label{restrictedammsequence} Let $d_j$ be the number of indecomposable components of $M^{\bf R}$ that correspond to the connected components of the $j$-th row of $Q(M^{\bf R})$. It is possible to define an action of an algebraic torus $T'=(\C^*)^{d_0-d_j}$ over $\mrm{Gr}_{\bf{e}}(M^{\bf R'})$  restricting the $T$ action over $\mrm{Gr}_{\bf{e}}(M^{\bf R})$. The $T'$-action endows ${\mathcal{F}l}^{\bf R'}(\C^n)$ with the structure of GKM variety and the $\bf{R'}$-admissible sequence $S'=(S'_1, \dots, S'_{n-1})$ then corresponds to a $T'$-fixed point $p_{S'}$.
\end{rmk}
\begin{rmk}
By Theorem~\ref{thm:flat}, if the linear degeneration ${\mathcal{F}l}^{\bf R}(\C^{n+1})$ is flat, the degeneration ${\mathcal{F}l}^{\bf R'}(\C^n)$ is also flat.
\end{rmk}
It is natural to ask how the mutations of $Q_S$ are linked to the ones of $Q_{S'}$. The following lemma explains this relation, crucial to prove our results about mutation graphs.
\begin{lemma}\label{restrictedmut}
A mutation from the $h$-th row of $Q_{S'}$ to the $k$-th one corresponds to a unique mutation from  $h$-th row to the $k$-th one of $Q_{S}$.
\end{lemma}
\begin{proof}
First of all, we remark that because of our construction, we can identify the rows of $Q_{S'}$ with the corresponding rows of $Q_S$. In particular, supposing that $Q_{S'}$ is obtained from $Q_S$ deleting the $j$-th row and the last column, if we have a mutation $[a,b]_{S'}^h \rightarrow k$ then $k$ must be different from $j$. Observe now that a mutation of the form $[a,b]_{S'}^h \rightarrow k$ with $b < n-1$ corresponds to a unique mutation of $Q_S$ because by definition of $R'$, if $b<n$ the segment $[a,b]_{M^{\bf R'}}^h$ is connected if and only if $[a,b]_{M^{\bf R}}^h$ is connected and $[a,b]_{S'}^h$ is a movable part if and only if $[a,b]_{S}^h$ is a movable part. Moreover, by Remark \ref{segmentmutations}, a mutation $[a,n-1]_{S'}^h \rightarrow k$ lift to a mutation $[a,n-1]_{S}^h \rightarrow k$ if and only if $\pi_{R_{n-1}}(k)=0$ or $k \in S_n$. The latter condition is always satisfied because $|S_n|=n$ and the Lemma follows.
\end{proof}
A consequence of the previous Lemma is that $\widetilde{G}_{S'}$ embeds into $\widetilde{G}_{S}$ by the restriction process, i.e. there exists a subgraph of $\widetilde{G}_{S}$, spanned by vertices $\{1, \dots, n+1\}\setminus \{j\}$  that is isomorphic to $\widetilde{G}_{S'}$.
As a consequence, if $\widetilde{G}_{S}$ is a complete graph, then $\widetilde{G}_{S'}$ must be again a complete graph, because $|\mut (p_{S'})|$ is always greater or equal than $\frac{n(n-1)}{2}$.\\
\begin{rmk}
A closer analysis of the mutations of $Q_S$ and $Q_{S'}$ implies that $\widetilde{G}_{S'}$ is a full subgraph of $\widetilde{G}_S$. We omit the proof of this fact since it is not necessary to prove our results. \end{rmk}
\begin{example}\label{ex:nonfull}
On the left of Figure~\ref{fig:nonfull} are displayed a fixed point $p_S$ for the $T$-action on the flat degeneration ${\mathcal{F}l}^{\bf R}(\C^{5})$ with ${\bf{R}}=(\{3,4\},\{3\},\{2,5\})$. On the right there is the successor closed subquiver corresponding to the fixed point $p_{S'}$ obtained from $p_S$ by restriction process. In particular, as observed in the Remark~\ref{restrictedammsequence} it is a fixed point for the $T'$ action on ${\mathcal{F}l}^{\bf R'}(\C^{5})$ where ${\bf{R}'}=(\{3,4\},\{3\})$. In Figure~\ref{fig:nonfullgraph} the oriented mutation graphs of $p_S$ and of $p_{S'}$ are displayed. In particular the blue full subgraph of $\widetilde{G}_S$ corresponds to the immersion of the oriented mutation graph $\widetilde{G}_{S'}$ induced by the restriction process. 
\end{example}
\begin{figure}[ht]
\centering
\begin{subfigure}{.3\textwidth}
\centering
\begin{tikzpicture}[baseline=(current bounding box.center), scale =0.9]
			\tikzstyle{point}=[circle,thick,draw=black,fill=black,inner sep=0pt,minimum width=2pt,minimum height=2pt]

			\node[] (v00) at (-1,1) {};
            \node[] (v01) at (-1,0.5) {};
			\node[] (v02) at (-1,0) {};
			\node[] (v03) at (-1,-0.5) {};
			\node[] (v04) at (-1,-1) {};
            \draw[fill] (v00)  circle (.05);
			\draw[fill] (v01)  circle (.05);
			\draw[fill] (v02)  circle (.05);
			\draw[fill,Red] (v03)  circle (.05);
			\draw[fill] (v04)  circle (.05);
			\node[] (v10) at (0,1) {};
            \node[] (v11) at (0,0.5) {};
			\node[] (v12) at (0,0) {};
			\node[] (v13) at (0,-0.5) {};
			\node[] (v14) at (0,-1) {};
            \draw[fill] (v10)  circle (.05);
			\draw[fill] (v11)  circle (.05);
			\draw[fill, Red] (v12)  circle (.05);
			\draw[fill,Red] (v13)  circle (.05);
			\draw[fill] (v14)  circle (.05);
			\node[] (v20) at (1,1) {};
            \node[] (v21) at (1,0.5) {};
			\node[] (v22) at (1,0) {};
			\node[] (v23) at (1,-0.5) {};
			\node[] (v24) at (1,-1) {};
            \draw[fill] (v20)  circle (.05);
			\draw[fill,Red] (v21)  circle (.05);
			\draw[fill,Red] (v22)  circle (.05);
			\draw[fill,Red] (v23)  circle (.05);
			\draw[fill] (v24)  circle (.05);
			\node[] (v30) at (2,1) {};
            \node[] (v31) at (2,0.5) {};
			\node[] (v32) at (2,0) {};
			\node[] (v33) at (2,-0.5) {};
			\node[] (v34) at (2,-1) {};
            \draw[fill] (v30)  circle (.05);
			\draw[fill,Red] (v31)  circle (.05);
			\draw[fill,Red] (v32)  circle (.05);
			\draw[fill,Red] (v33)  circle (.05);
			\draw[fill,Red] (v34)  circle (.05);
			\draw[, ] (v00) -- (v10);
			\draw[thick, Red ] (v03) -- (v13);
			\draw[, ] (v04) -- (v14);
			\draw[, ] (v10) -- (v20);
			\draw[,  ] (v11) -- (v21);
			\draw[thick,Red ] (v13) -- (v23);
			\draw[, ] (v14) -- (v24);
			\draw[thick,Red ] (v22) -- (v32);
			\draw[thick,Red  ] (v21) -- (v31);
			\draw[,] (v24) -- (v34);
				\draw[,] (v21) -- (v31);
		\end{tikzpicture}
\end{subfigure}\hspace{0.1\textwidth}
\begin{subfigure}{.3\textwidth}
\centering
\begin{tikzpicture}[baseline=(current bounding box.center), scale =1.2]
			\tikzstyle{point}=[circle,thick,draw=black,fill=black,inner sep=0pt,minimum width=2pt,minimum height=2pt]

            \node[] (v01) at (-1,0.5) {};
			\node[] (v02) at (-1,0) {};
			\node[] (v03) at (-1,-0.5) {};
			\node[] (v04) at (-1,-1) {};
			\draw[fill] (v01)  circle (.05);
			\draw[fill] (v02)  circle (.05);
			\draw[fill,Red] (v03)  circle (.05);
			\draw[fill] (v04)  circle (.05);
            \node[] (v11) at (0,0.5) {};
			\node[] (v12) at (0,0) {};
			\node[] (v13) at (0,-0.5) {};
			\node[] (v14) at (0,-1) {};
			\draw[fill] (v11)  circle (.05);
			\draw[fill, Red] (v12)  circle (.05);
			\draw[fill,Red] (v13)  circle (.05);
			\draw[fill] (v14)  circle (.05);
            \node[] (v21) at (1,0.5) {};
			\node[] (v22) at (1,0) {};
			\node[] (v23) at (1,-0.5) {};
			\node[] (v24) at (1,-1) {};
			\draw[fill,Red] (v21)  circle (.05);
			\draw[fill,Red] (v22)  circle (.05);
			\draw[fill,Red] (v23)  circle (.05);
			\draw[fill] (v24)  circle (.05);
			\draw[, ] (v04) -- (v14);
			\draw[,  ] (v11) -- (v21);
			\draw[thick,Red ] (v13) -- (v23);
           \draw[thick,Red  ] (v03) -- (v13);
           	\draw[thick,Red ] (v13) -- (v23);
			\draw[,] (v14) -- (v24);
		\end{tikzpicture}
\end{subfigure}
\caption{A fixed point $p_S$ (on the left) and its restriction $p_{S'}$ (on the right).}\label{fig:nonfull}
\end{figure}
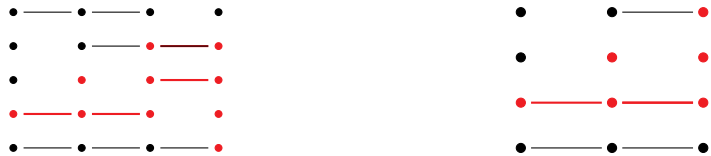
\begin{figure}[ht]
\centering
\begin{subfigure}{.3\textwidth}
\centering
\begin{tikzpicture}
\node (a) at (0,0) {};
\node (b) at (4,0) {};
\node (c) at (2,1){};
\node (d) at (2,3){};
\node (e) at (0,2){};

\draw[fill=black] (0,0) circle (2pt);
\draw[fill=black] (4,0) circle (2pt);
\draw[fill=black] (2,1) circle (2pt);
\draw[fill=black] (2,3) circle (2pt);
\draw[fill=black] (0,2) circle (2pt);
\node at (-0.3,0) {4};
\node at (4.3,0) {3};
\node at (2.3,1.2) {2};
\node at (2,3.3) {1};
\node at (-0.3,2.3) {5};
    \draw[thick,blue,-latex] (b) to[bend left]  (a);
    \draw[thick,blue,-latex](c) -- (a);
    \draw[thick,blue,-latex] (a) -- (d);
     \draw[thick,blue,-latex] (c) -- (b);
    \draw[thick,blue,-latex](b) to  (d);
    \draw[thick,blue, -latex] (c) -- (d);
\draw [thick, Red, -latex] (b) to [bend left] (e);
        \draw[thick,Red,-latex](d) to (e);
            \draw[thick,Red,-latex](c) to (e);
                \draw[thick,Red,-latex](a) to (e);
\end{tikzpicture}
\end{subfigure}\hspace{0.05\textwidth}
\begin{subfigure}{.3\textwidth}
\centering
\begin{tikzpicture}
\node (a) at (0,0) {};
\node (b) at (4,0) {};
\node (c) at (2,1){};
\node (d) at (2,3){};

\draw[fill=black] (0,0) circle (2pt);
\draw[fill=black] (4,0) circle (2pt);
\draw[fill=black] (2,1) circle (2pt);
\draw[fill=black] (2,3) circle (2pt);
\node at (-0.3,0) {4};
\node at (4.3,0) {3};
\node at (2.3,1.2) {2};
\node at (2,3.3) {1};
    \draw[thick,-latex] (b) to  (a);
    \draw[thick,-latex](c) -- (a);
    \draw[thick,-latex] (a) -- (d);
    \draw[thick,-latex](b) to (d);
    \draw[thick,-latex] (c) -- (d);
\draw [thick, -latex] (c) to  (b);
\end{tikzpicture}
\end{subfigure}
\caption{The oriented mutation graphs of $p_S$ and of $p_{S'}$}\label{fig:nonfullgraph}
\end{figure}
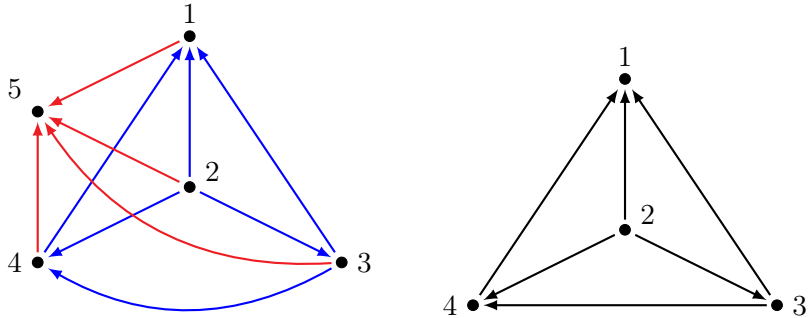
\subsection{Smoothness and Mutation Graph}\label{subsec:proofgraph}
In this section we prove that if a point $p_S$ is smooth, then the oriented mutation graph $\widetilde{G}_S$ is a transitive tournament.
We deduce preliminary informations about the orientation of $\widetilde{G}_S$ as a consequences of Lemmata~\ref{emptyrow} and \ref{emptymut} and we use the tools developed in Section~{\ref{Restriction}} to show that the condition $(1)$ in Theorem~{\ref{thm:criteria}} implies that the condition $(5)$ holds. \\
Firstly we remark that the mutation graph $G_S$ is the complete graph over $n+1$ vertices if and only if for every pair of indices $1 \leq i,j \leq n+1$
\[\mathrm{Mut}_S(i,j) + \mathrm{Mut}_S(j,i) =1 \]
To prove that conditions $(1) $ and $(5)$ in Theorem~\ref{thm:criteria} are equivalent, we need two preliminary Lemmata:
\begin{lemma}\label{emptyrow}
If $p_S$ is smooth, then $Q_S$ has an empty row.
\end{lemma}
\proof
Suppose that $Q_S$ does not have empty rows. By Remark \ref{HomTg}, it is enough to prove that under this assumption $\Ext^1(M_S, M^{\bf R}/M_S)\neq 0$.  By definition of admissible sequence and because $Q_S$ does not have an empty row, there exists $i \in \{1, \dots, n+1\}\setminus S_n$ and $h$  such that $i \in S_h$  and $i \notin S_j$ for all $j \geq h$. Set $k_i = \min\{ k \in \{1, \dots, n\} \, | \, i \in S_j,  \; \forall j \; k \leq j \leq h\} $ and $m_i=\max \{m \in \pos(i) \, | \, m < h\}$, the segment $[\max \{k_i,m_i\}, h]^i_S$ is a connected component of $Q_S$ ad it corresponds to the indecomposable summand $U^i_{\max \{k_i,m_i\}, h}$  of $M_S$. Moreover, if $P=\min \{p \in \pos(i) \cup \{n\} \, | \, h < p\}$, we have that $[h+1, P]^i_S$ is a connected component of $Q(M^{\bf R})\setminus Q_S$ and the representation $U_{h+1, P}$ is an indecomposable summmand of $M^{\bf R}/M_S$. By Formula~(\ref{ext}) we have $\Ext^1(U_{\max \{k_i,m_i\}, h}, U_{h+1, P})\neq 0$ and $p_S$ cannot be smooth.
\endproof
\begin{lemma}\label{emptymut}Let $p_S$ be a fixed point for the degeneration ${\mathcal{F}l}^\mathbf{R}(\C^{n+1})$. If the $j$-th row of $Q_S$ is empty, then $|\mut_S(i,j)|\geq 1$ for all $i \in \{1, \dots, n+1\}, \; i \neq j$. 
\end{lemma}
\begin{proof}
First of all observe that such a $j$ is unique, because $|S_n|=n$. For $i \neq j$ set 
$m_i=\max\{ \pos(i) \cup \{0 \}\}+1$
and denote by $k_i$ the minimum of the set $\{h | i \in S_h, m_i \leq h\}$ (observe that $k_i$ exists because $j$ is the unique index not in $ S_n$). Set now 
$m(i,j)= \min \{k \in \pos(j)\cup \{n+1\}\; | \; k \geq k_i\}$. By definition  $[{k_i},{m(i,j)}]_S^i$ is a movable part and $[{k_i},{m(i,j)}]_{M^{\bf R}}^j$ is a connected segment with ending point equal to $n+1$ or to an element in $\pos(j)$. 
This implies, by Remark \ref{segmentmutations}, that we have a mutation $[{k_i},{m(i,j)}]_S^i\rightarrow j$ and the lemma is proved.
\end{proof}
\begin{rmk}\label{rmk:sink}
Previous Lemmata highlight also some properties of the oriented mutation graph.
In fact Lemma~\ref{emptymut} implies that if the $j$-th row of $Q_S$ is empty, then the vertex $j$ is a sink for $\widetilde{G}_S$.
\end{rmk}
Moreover, in our proof we need some classical results about transitive tournaments. We refer to \cite{HM} for a complete survey on the Tournament Graphs and related properties.
\begin{defn}
Let $v$ be a vertex of a directed graph $G$, the indegree $\mathrm{id}_G(v)$ (rep. outdegree $\mathrm{od}_G(v)$) of $v$ is the number of edges $e \in E(G)$ such that $v$ is the target (resp. the source) of $e$. 
\end{defn}
The property of being transitive is linked to acyclicity and to the sets of outdegrees and indegrees, respectively.
\begin{theorem}[\cite{HM}, Corollary 5a]\label{TT}
Let $T$ be a tournament graph over $n+1$ vertices. The following are equivalent:
\begin{itemize}
    \item $T$ is a transitive tournament,
    \item $T$ is acyclic,
     \item  $T$ does not contain a cycle of length 3,
    \item The set of outdegrees of $T$ is $\{0, \dots, n\}$.
     \item The set of indegrees of $T$ is $\{0, \dots, n\}$.
\end{itemize}
\end{theorem}
In particular a full subgraph of a transitive tournament is again a transitive tournament.
\begin{proposition}
If the fixed point $p_S$ for ${\mathcal{F}l}^{\bf{R}}(\C^n)$ is smooth then $\widetilde{G}_S$ is a transitive tournament over $n$ vertices.
\end{proposition}
\proof
We proceed by induction on $n$.
The base case $n=1$ is completely trivial. For the general case, firstly we remark that by Lemma~\ref{emptyrow}, there exists $j$ such that $j \notin S_h$ for all $h\leq n$. Using the restriction process described in Section~\ref{Restriction}, we can consider an $\bf{R'}$-admissible sequence $S'$ for the degeneration associated to the quiver $Q(M^{\bf R'})$ obtained by $Q(M^{\bf R})$ deleting the $j$-th row and the last column.
We can easily obtain the that $p_{S'}$ is smooth: as a consequence of Lemma~\ref{restrictedmut} we have \[ |\mut(p_{S'})| + \sum_{i \neq j} |\mut_S(i,j)| \leq \sum_{h,k \neq j} |\mut_S(h,k)|+ \sum_{i \neq j} |\mut_S(i,j)| = \frac{n(n+1)}{2},\]
moreover, Lemma \ref{emptymut} implies that $\sum_{i \neq j} |\mut_S(i,j)|\geq n$ and consequently 
 \[\frac{n(n-1)}{2} \leq |\mut(p_{S'})| \leq \frac{n(n+1)}{2} - n.\]
So, by inductive reasoning, we can suppose $\widetilde{G}_{S'}$ to be a transitive tournament over $n$ vertices. The previous computation, combined with the Lemma \ref{emptymut}, implies that $|\mut_S(i,j)|=1$ for every $i \neq j$. As an immediate consequence, the mutation graph $G_S$ is a complete graph. Moreover structure of $\widetilde{G}_S$ is more clear: its full subgraph spanned by the vertices $V(\widetilde{G}_S) \setminus \{j\}$ can be identified with $\widetilde{G}_{S'}$ and the vertex corresponding to $j$ is connected to each other by a single edge. In particular, by Remark~{\ref{rmk:sink}}, the vertex $j$ is a sink for $\widetilde{G}_S$ and consequently $\mrm{id}_{\widetilde{G}_S}(j)=n$. Moreover, the fact that $j$ is a sink implies that the indegree of a vertex $v\neq j$ in $\widetilde{G}_S$ is equal to $\mrm{id}_{\widetilde{G}_{S'}}(v)$. Consequently the set of indegrees of $\widetilde{G}_S$ is $\{0,\dots, n\}$ and by Theorem~\ref{TT} we obtain that $\widetilde{G}_S$ is a transitive tournament over $n+1$ vertices. 
\endproof

\end{document}